\documentclass[12pt]{article}
\usepackage{mathrsfs}
\usepackage{amssymb} \textwidth 150mm \textheight 220mm \oddsidemargin
15pt \evensidemargin 0pt \topmargin 0cm \headsep 0.3cm

\usepackage{amsmath}
\usepackage{graphicx}
\usepackage{indentfirst}
\usepackage{extarrows}
\usepackage{cite}
\usepackage{float}
\usepackage{dsfont}
\usepackage{CJK}
\usepackage{multirow}
\usepackage{makecell}

\usepackage{subfigure}

\begin{document}
\title{\Large\bf{Zeros of complete elliptic integrals and its application to Melnikov functions}}
\author{{Jihua Yang\thanks{E-mail addresses: jihua1113@163.com, yangjh@mail.bnu.edu.cn}}\\
{\small \it School of Mathematical Sciences, Tianjin Normal University, Tianjin 300387, China}\\ }
\date{}
\maketitle \baselineskip=0.9\normalbaselineskip \vspace{-3pt}

\noindent
{\bf Abstract}\, In this paper, we first discuss the linear independence of the complete elliptic integrals of the first, second and third kinds $K(k)$, $E(k)$ and $\Pi(\mu(k),k)$, and then obtain an upper bound for the number of zeros of a function of the form
\begin{eqnarray*}
p(k)K(k)+q(k)E(k)+r(k)\Pi(\mu(k),k),\ k\in(-1,1),
\end{eqnarray*}
where $p(k)$, $q(k)$ and $r(k)$ are real polynomials, $\mu(k)$ is a real polynomial or rational function. Finally, we apply it to a Hamiltonian triangle with three invariant straight lines under small real polynomials piecewise smooth perturbation.

\vskip 0.2 true cm


 \section{Introduction and main results}
 \setcounter{equation}{0}
\renewcommand\theequation{1.\arabic{equation}}

Consider the perturbed Hamiltonian system
\begin{eqnarray}
\begin{cases}
\dot{x}  =  H_y(x,y)+\varepsilon f(x,y),\\
\dot{y}  =  -H_x(x,y)+\varepsilon g(x,y),\\
\end{cases}
\end{eqnarray}
with $0<|\varepsilon|\ll1$, where $H(x,y)$ is a real  polynomial  of degree $m+1$, and $f(x,y)$ and $g(x,y)$ are real  polynomials of degree $n$ in the variables $x$ and $y$. When $\varepsilon=0$,  system (1.1) is assumed to admit a family of closed orbits $\{\Gamma_h\}$, parameterized by $h\in\Sigma$, with $\Sigma$ as the maximal  open  interval where $h$ exists. That is,
$$\Gamma_h=\{(x,y)\in\mathds{R}^2\,|\,H(x,y)=h,\,h\in\Sigma\}.$$
The first order Melnikov function $I(h)$ of system (1.1) was obtained by Pontryagin \cite{PL} which is given by the Abelian integral
\begin{eqnarray*}
I(h)=\oint_{\Gamma_h}g(x,y)dx-f(x,y)dy.
\end{eqnarray*}
For given positive integers $m,n\geq2$, %all $H(x,y)$, $f(x,y)$ and $g(x,y)$ satisfying the above conditions,
and all possible closed orbits $\Gamma_h$, the problem of determining the maximum number of isolated zeros of the first order Melnikov function $I(h)$ (when it is not identically zero) is known as the weak Hilbert's 16th problem \cite{A,AV}.  After decades of relentless efforts by scholars, the weak Hilbert's 16th problem has been completely solved for the case where $m=n=2$, see \cite{M,ZL,LZ,G,HI,HI94,CLLZ}. However, when $m\geq3$, the problem is far from being resolved, although many significant results have been achieved, see \cite{BNY,HI98,YZ17,ZZ,CZW,XH20,DL1,DL2,DL3,DL4,GGM16,MV,CH,LX,LX13} and the references therein.
%Considerable research efforts have been devoted to the weak Hilbert's 16th problem. Notable works in this area include, but are not limited to, \cite{BNY,HI98,YZ,YZ17,ZZ,CZW,NM,XH20,DL1,DL2,DL3,DL4,GGM16}.

When system (1.1) has a switching straight line $y=cx+d$, $c,d\in\mathds{R}$, and
 $$f(x,y)=\begin{cases}
 f^+(x,y),\ y\geq cx+d,\\
 f^-(x,y),\ y< cx+d,
 \end{cases}
 g(x,y)=\begin{cases}
 g^+(x,y),\ y\geq cx+d,\\
 g^-(x,y),\ y< cx+d,
 \end{cases}
 $$
 where $f^\pm(x,y)$ and $g^\pm(x,y)$ are real polynomials in $x$ and $y$ of degree $n$. The corresponding first order Melnikov function $I(h)$ was obtained by Han et al. \cite{LH,LH12} as follows
\begin{eqnarray}
I(h)=\int_{\Gamma^+_h}g^+(x,y)dx-f^+(x,y)dy+\int_{\Gamma^-_h}g^-(x,y)dx-f^-(x,y)dy,
\end{eqnarray}
where
$$\begin{aligned}
&\Gamma^+_h=\{(x,y)\in\mathds{R}^2\,|\,H(x,y)=h,\,h\in\Sigma,\,y\geq ax+b\},\\
&\Gamma^-_h=\{(x,y)\in\mathds{R}^2\,|\,H(x,y)=h,\,h\in\Sigma,\,y< ax+b\}.
\end{aligned}$$
There has also been a wealth of significant research on estimating the number of zeros of Melnikov functions for perturbed differential systems with switching lines; see, for instance, \cite{YZ,LL,XH,ZSJ,BBLN,LHR,TH,LHL,PR}.

 The relationship between the number of limit cycles in system (1.1) and the number of isolated zeros of the corresponding first order Melnikov function is well known \cite{G,IY}. Specifically, the total number of zeros (counted with multiplicity) of the first order Melnikov function $I(h)$ constitutes an upper bound for the limit cycles bifurcating from an open period annulus $\bigcup\limits_{h\in\Sigma}\Gamma_h$ or $\bigcup\limits_{h\in\Sigma}\Gamma^+_h\cup\Gamma^-_h$ \cite{G,HS}. This result extends to a closed period annulus bounded by a homoclinic loop, as proven by Roussarie \cite{R}. Furthermore, the number of simple multiple zeros of $I(h)$ gives a lower bound for the number of limit cycles.

The calculation of a first order Melnikov function usually  involves elliptic integrals, which can be expressed in terms of the complete elliptic integrals of the first, second, and third kinds:
\begin{eqnarray}\begin{aligned}
&K(k)=\int^1_0\frac{1}{\sqrt{(1-z^2)(1-k^2z^2)}}dz,\\
&E(k)=\int^1_0\frac{\sqrt{1-k^2z^2}}{\sqrt{1-z^2}}dz,\\
&\Pi(\mu(k),k)=\int^1_0\frac{1}{(1-\mu z^2)\sqrt{(1-z^2)(1-k^2z^2)}}dz.
\end{aligned}\end{eqnarray}
Here, $k$ is referred to as the modulus with its value typically confined to the interval $(-1,1)$, and $\mu(k)$ is called the parameter of the integral of the third kind. %Moreover, $\mu(k)$ is a nonnegative real function of $k$, and $\mu(k)\neq k^2,\mu(k)\neq1.$
Therefore, it is important to estimate the number of zeros for a function of the form
\begin{eqnarray}
I(k)=p(k)K(k)+q(k)E(k)+r(k)\Pi(\mu(k),k),\ k\in(-1,1),
\end{eqnarray}
where $p(k)$, $q(k)$ and $r(k)$ are real polynomials of $k$ and $\mu(k)$ is a real polynomial or rational function of $k$.
When $r(k)=0$, Gasull, Li, Llibre and Zhang proved an important and well-known result in \cite{GLLZ} by using the Argument Principle. For convenience, we state it as the theorem below.

 \vskip 0.2 true cm

\noindent
{\bf Theorem 1.1}\, \cite{GLLZ} {\it  Let $p(k)$ and $q(k)$ be real polynomials in $k$ of degrees $m$ and $n$, respectively and $k\in(-1,1)$. Then an upper bound for the number of zeros of the function
\begin{eqnarray}I(k)=p(k)K(k)+q(k)E(k),\end{eqnarray}
taking into
account their multiplicities, is $m+n+2$. } \vskip 0.2 true cm

In the present paper, we deal with the case where $r(k)\neq0$. This constitutes the primary objective of the first part of this paper. The desired estimate relies on establishing some properties for $K(k)$, $E(k)$ and $\Pi(\mu(k),k)$ \cite{C,CG,GMV,GGM,GI,GI03,GI09}. The main results are presented in the following four theorems.
\vskip 0.2 true cm

\noindent
{\bf Theorem 1.2}\, {\it Let $p(k)$, $q(k)$, $r(k)$ and $\mu(k)$  be differentiable real functions.  Then there exist real functions $M(k)$ and $N(k)$ such that
\begin{eqnarray*}r^2(k)\frac{d}{dk}\Big(\frac{I(k)X(k)}{r(k)}\Big)=M(k)K(k)+N(k)E(k),\ k\in(-1,1)\backslash S,\end{eqnarray*}
where $S$ is the set  of  zeros of $r(k)$, and $X(k)$ is defined by (2.12).}
\vskip 0.2 true cm

According to Theorem 1.2, if an upper bound for the number of zeros of $M(k)K(k)+N(k)E(k)$ can be established on the interval $(-1,1)$, then  an upper bound for the number of zeros of $I(h)$ on $(-1,1)$  follows via Rolle's Theorem. However, this is a nontrivial task for general functions $M(k)$ and $N(k)$. Nevertheless, it can be achieved when $p(k)$, $q(k)$, $r(k)$ and  $\mu(k)$ are  real polynomials or rational functions. This is because the properties of $M(k)$ and $N(k)$  are fully determined by those of $p(k)$, $q(k)$, $r(k)$ and  $\mu(k)$.
 \vskip 0.2 true cm

\noindent
{\bf Theorem 1.3}\, {\it  Let $p(k)$, $q(k)$, $r(k)$ and $\mu(k)$ be real polynomials in $k$ of degrees $m$, $n$, $l$ and $s$, respectively and $k\in(-1,1)$. Then the number of zeros of $I(k)$ in (1.4) does not exceed   $\psi(m,n,l,s)$, taking into account their multiplicities, where
\begin{eqnarray*}
\psi(m,n,l,s)= \begin{cases}\begin{cases}
2\max\{m,n\}+3l+6s+7,\ \textup{if}\ s\geq2,\,l\leq \max\{m,n\}+s,\\
5l+4s+7,\qquad\qquad\qquad\ \ \, \textup{if}\ s\geq2,\,l>\max\{m,n\}+s,\\
\end{cases}\\
\begin{cases}
2\max\{m,n\}+3l+15,\qquad \textup{if}\ s=1,\, l\leq \max\{m,n\}+1,\\
5l+13,\qquad\qquad\qquad\qquad\ \, \textup{if} \ s=1,\, l>\max\{m,n\}+1,\\
\end{cases}\\
\begin{cases}
2\max\{m,n\}+3l+11,\qquad \textup{if}\ s=0,\, l\leq \max\{m,n\}+2,\\
\max\{m,n\}+4l+9,\quad\quad\ \ \ \, \textup{if} \ s=0,\, l>\max\{m,n\}+2.
\end{cases}\\
\end{cases}
\end{eqnarray*}
Moreover, there exist real polynomials $p(k)$, $q(k)$ and $r(k)$ such that $I(k)\not\equiv0$ and  $I(k)$ has $m+n+l+2$ zeros on the interval $(-1,1)$.}
 \vskip 0.2 true cm

 When $p(k)$, $q(k)$ and $r(k)$  are  real polynomials  and   $\mu(k)$ is a  real rational function, the upper bound of the number of zeros of $I(k)$ is also can be obtained.

 \vskip 0.2 true cm

\noindent
{\bf Theorem 1.4}\, {\it  Let $p(k)$, $q(k)$ and $r(k)$   be real polynomials of $k$ of degrees $m$, $n$ and $l$, and $\mu(k)=\frac{2k^2}{1+k^2}$. Then the number of zeros of $I(k)$ in (1.4) does not exceed $\psi(m,n,l)$, taking into account their multiplicities, where
\begin{eqnarray*}
\psi(m,n,l)= \begin{cases}
2\max\{m,n\}+3l+11,\ \textup{if}\ l\leq \max\{m,n\}+2,\\
5l+7,\qquad\qquad\qquad\quad \,\, \textup{if}\ l>\max\{m,n\}+2.\\
\end{cases}
\end{eqnarray*}Moreover, there exist real polynomials $p(k)$, $q(k)$ and $r(k)$ such that $I(k)\not\equiv0$ and  $I(k)$ has $m+n+l+2$ zeros on the interval $(-1,1)$.}

 \vskip 0.2 true cm

\noindent
{\bf Remark 1.1}\, In fact, when $\mu(k)=\frac{Q(k)}{P(k)}$, where $P(k)$ and $Q(k)$ are real polynomial of $k$, the upper bound for the number of zeros of $I(k)$ on the interval $(-1,1)$ can also be obtained, similarly. It is just that the calculation process and results are more complex.

 \vskip 0.2 true cm

As an application of Theorem 1.4, we study the number of zeros of the first order Melnikov function $I(h)$ associated to the following piecewise smooth perturbed Hamiltonian triangle
\begin{eqnarray}
\left(  \begin{array}{c}
          \dot{x} \\          \dot{y}
          \end{array} \right)
=\begin{cases}
 \left(
  \begin{array}{c}
          x^2(1-x-2y)+\varepsilon f^+(x,y) \\
          xy(-2+3x+2y)+\varepsilon g^+(x,y)
          \end{array} \right), \quad y\geq\frac12-\frac12x,\\[0.6truecm]
 \left(  \begin{array}{c}
         x^2(1-x-2y)+\varepsilon f^-(x,y) \\
          xy(-2+3x+2y)+\varepsilon g^-(x,y)
           \end{array}
 \right),\quad y<\frac12-\frac12x,
  \end{cases}
 \end{eqnarray}
 where $0<|\varepsilon|\ll1$,
$$\begin{aligned}
&f^\pm(x,y)=\sum\limits_{i+j=0}^na^\pm_{i,j}x^iy^j,\ g^\pm(x,y)=\sum\limits_{i+j=0}^nb^\pm_{i,j}x^iy^j,\ a^\pm,b^\pm\in\mathds{R}, i,j\in\mathds{N}.
\end{aligned}$$

\noindent
{\bf Theorem 1.5}\, {\it  An upper bound for the number of zeros of the first order Melnikov function associated to system (1.6) is  $\big[\frac{11n}{2}\big]+43$,  taking into account their multiplicities.}

\vskip 0.2 true cm

This paper is organized as follows. In Section 2, we first present the differential equations satisfied by $K(k)$, $E(k)$ and $\Pi(\mu(k),k)$, as well as their linear independence, and then prove Theorems 1.2-1.4. The proof of  Theorem 1.5 is given in Section 3. Finally, a conclusion is drawn in Section 4.

 \section{Elliptic integrals}
 \setcounter{equation}{0}
\renewcommand\theequation{2.\arabic{equation}}

This paper considers the problem in the real domain. According to the definition of the complete elliptic integral of the third kind, it can be seen that $\mu(k)$ is a non-negative real function of $k$, with $\mu(k)\neq k^2$ and $\mu(k)\neq1$. It should be noted that $\mu(k)$ may or may not be a function of $k$. When $\mu(k)=1$, the integral $\Pi(\mu(k),k)$ is divergent.  When $\mu(k)=k^2$, it is straightforward  to check that
\begin{eqnarray}
\Pi(\mu(k),k)=\frac{1}{1-k^2}E(k).
\end{eqnarray}
At this point, the first Melnikov function $I(k)$ in (1.4) can be transformed into the form of (1.5), and the corresponding estimate of the number of zeros can be obtained immediately. Hence, this case will not be addressed in the present work.

 \subsection{Properties of elliptic integrals}

 We first give the differential equations associated with the complete elliptic integrals $K(k)$, $E(k)$ and $\Pi(\mu(k),k)$. They play  a very crucial role in the proof of the main results.
\vskip 0.2 true cm

\noindent
 {\bf Lemma 2.1} {\it Let $\mu(k)$ be a twice continuously differentiable real function of $k$ and $V(k)=\big(K(k),E(k),\Pi(\mu(k),k)\big)^T$. Then the following hold:}
 \vskip 0.2 true cm

\noindent
(i) {\it The vector function $V(k)$ satisfies the following Picard-Fuchs equation
 \begin{eqnarray}
 V'(k)=\mathbf{A}(k) V(k),\ k\in(-1,1),
\end{eqnarray}
where
 \begin{eqnarray*}
\mathbf{A}(k)=\left(\begin{matrix}
                -\frac{1}{k}&\frac{1}{k(1-k^2)}&0\\
                 -\frac{1}{k}&\frac{1}{k}&0\\
                 a_{31}(k)&a_{32}(k)&a_{33}(k)\\
\end{matrix}\right),
\end{eqnarray*}
\begin{eqnarray*}\begin{aligned}
&a_{31}(k)=-\frac{\mu'(k)}{2\mu(1-\mu)},\\ &a_{32}(k)=\frac{k}{(\mu-k^2)(k^2-1)}+\frac{\mu'(k)}{2(1-\mu)(\mu-k^2)},\\&
a_{33}(k)=\frac{k}{\mu-k^2}+\frac{(\mu^2-k^2)\mu'(k)}{2\mu(1-\mu)(\mu-k^2)}.
\end{aligned}\end{eqnarray*}
}
 \vskip 0.2 true cm

\noindent
(ii) {\it The vector function $V(k)$ satisfies the following linear differential equations of order two
\begin{eqnarray}
V''(k)=\mathbf{A}(k)V'(k)+\mathbf{B}(k)V(k),\ k\in(-1,1),
\end{eqnarray}
where
 \begin{eqnarray*}
\mathbf{B}(k)=\left(\begin{matrix}
                \frac{1}{k^2}&\frac{3k^2-1}{k^2(1-k^2)^2}&0\\
                 \frac{1}{k^2}&-\frac{1}{k^2}&0\\
                 b_{31}(k)&b_{32}(k)&b_{33}(k)\\
\end{matrix}\right),
\end{eqnarray*}
$$\begin{aligned}
b_{31}(k)=&\frac{(\mu(k)^2-\mu(k))\mu''(k)-(2\mu(k)-1)\mu'(k)^2}{2\mu(k)^2(\mu(k)-1)^2},\\
b_{32}(k)=&\frac{3k^4-\mu(k)(k^2+1)-k^2}{(1-k^2)^2(\mu(k)-k^2)^2}+\frac{1}{2(1-\mu(k))(\mu(k)-k^2)}\mu''(k)\\&
+\frac{2\mu(k)-k^2-1}{2(\mu(k)-1)^2(\mu(k)-k^2)^2}\mu'(k)^2-\frac{k(k^2+\mu(k)-2)}{(k^2-1)(\mu(k)-k^2)^2(\mu(k)-1)}\mu'(k),\\
b_{33}(k)=&\frac{\mu(k)+k^2}{(\mu(k)-k^2)^2}-\frac{k^2-\mu(k)^2}{2\mu(k)(\mu(k)-1)(k^2-\mu(k))}\mu''(k)\\&
+\frac{(2\mu(k)-1)k^4-2(2\mu(k)^2-\mu(k))k^2+\mu(k)^4}{2\mu(k)^2(\mu(k)-1)^2(k^2-\mu(k))^2}\mu'(k)^2-\frac{2k}{(\mu(k)-k^2)^2}\mu'(k).
\end{aligned}$$
}

\noindent
{\bf Proof}\,  (i) The first and second differential equations in (2.2) can be derived from formulas 710.00 and 710.02 in \cite{BF}. We only prove the third one in (2.2). For notational convenience, we denote $$F(\mu,k,z)=\frac{1}{(1-\mu z^2)\sqrt{(1-z^2)(1-k^2z^2)}}.$$
Evidently, $F(\mu,k,z)$ and $\frac{\partial F(\mu,k,z)}{\partial k}$ are continuous when $(k,z)\in (-1,1)\times(0,1)$ and $F(\mu,k,z)$ and $\frac{\partial F(\mu,k,z)}{\partial \mu}$ are continuous when $(\mu,z)\in (0,C]\times(0,1)$,  here $C$ is an arbitrary positive number. Hence,
\begin{eqnarray*}\begin{aligned}
&\frac{\partial \Pi(\mu,k)}{\partial k}=\int_0^1\frac{k z^2}{\left(1-k^2 z^2\right)\left(1 - \mu z^2\right)  \sqrt{\left(1 - z^2\right)\left(1 - k^2 z^2\right)}}dz,\\
&\frac{\partial \Pi(\mu,k)}{\partial \mu}=\int_0^1\frac{z^2}{\sqrt{\left(1-z^2 \right)\left(1-k^2 z^2 \right)} \left(1-\mu z^2 \right)^2}dz.
\end{aligned}\end{eqnarray*}
Observing that
$$\frac{z^2}{\left(1-k^2 z^2\right)\left(1-\mu z^2\right)}= \frac{1}{\left(k^2-\mu\right)\left(1-k^2 z^2\right)}+ \frac{1}{\left(\mu-k^2\right)\left(1-\mu z^2\right)},$$
one finds that
\begin{align*}
%I &= k \int_{0}^{1} \left[ \frac{1}{\left(k^2-a\right)\left(1-k^2 z^2\right)} + \frac{1}{\left(a-k^2\right)\left(1-a z^2\right)} \right] \frac{1}{\sqrt{\left(1-z^2\right)\left(1-k^2 z^2\right)}} dz \\
\frac{\partial \Pi(\mu,k)}{\partial k}= \frac{k}{k^2-\mu} \int_{0}^{1} \frac{1}{\left(1 - k^2 z^2\right)^{3/2} \sqrt{1 - z^2}} dz + \frac{k}{\mu-k^2} \Pi(\mu, k).
\end{align*}
It follows from (2.1) that
\begin{eqnarray}
\frac{\partial \Pi(\mu,k)}{\partial k}= \frac{k E(k)}{\left(1 - k^2\right)(k^2-\mu)} + \frac{k\Pi(\mu, k)}{\mu-k^2}.
\end{eqnarray}

On the other hand, a straightforward computation gives
\begin{eqnarray}
\frac{\partial \Pi(\mu,k)}{\partial \mu}=-\frac{1}{\mu}\Pi(\mu,k)+\frac1\mu \int_0^{\frac\pi2}\frac{1}{(1-\mu\sin^2\theta)^2\sqrt{1-k^2\sin^2\theta}}d\theta.
\end{eqnarray}
In light of the formula 336.02 in \cite{BF}, one finds that
$$\begin{aligned}\int_0^{\frac\pi2}\frac{1}{(1-\mu\sin^2\theta)^2\sqrt{1-k^2\sin^2\theta}}d\theta=&\frac{1}{2(\mu-1)(k^2-\mu)}\big[(k^2-\mu)K(k)+\mu E(k)\\&+(2\mu k^2+2\mu-\mu^2-3k^2)\Pi(\mu,k)\big].\end{aligned}$$
Therefore, inserting the above formula into (2.5) yields
\begin{eqnarray}\begin{aligned}
\frac{\partial \Pi(\mu,k)}{\partial \mu}= \frac{ K(k)}{2 \, \mu \, (\mu - 1)} - \frac{E(k)}{2 \, (\mu - 1) \, \left(\mu-k^2 \right)} - \frac{\left(\mu^2 - k^2\right) \, \Pi(\mu, \, k)}{2 \,  \mu \, (\mu - 1)\left(\mu-k^2 \right)}.
\end{aligned}\end{eqnarray}
By the chain rule for the differentiation of composite functions, it follows that
\begin{eqnarray}
\frac{d \Pi(\mu(k),k)}{d k}= \frac{\partial \Pi(\mu,k)}{\partial k}+\frac{\partial \Pi(\mu,k)}{\partial \mu}\frac{d \mu(k)}{d k}.
\end{eqnarray}
Substituting (2.4) and (2.6)  into (2.7) yields the desired result.

(ii) Differentiating both sides of (2.2) with respect to $k$ yields equation (2.3). This completes the proof. $\Box$
\vskip 0.2 true cm

In order to prove the main results of this paper, we also need to discuss the linear independence of $K(k)$, $E(k)$ and $\Pi(\mu(k),k)$. To this end, we introduce the symmetric integrals of the first, second and third kinds, see \cite{C87,C88,C89,C91,C92}. %Assume that $x$, $y$ and $z$ are nonnegative and at most one of them is $0$.
 The symmetric integrals of the first and third kinds are
$$R_F(x,y,z)=\frac12\int_0^\infty\frac{1}{\sqrt{(t+x)(t+y)(t+z)}}dt,$$
and
$$R_J(x,y,z,p)=\frac32\int_0^\infty\frac{1}{\sqrt{(t+x)(t+y)(t+z)}(t+p)}dt, \ p>0.$$
If $p=z$, then $R_J(x,y,z,p)$ reduces to the symmetric integral of the second kind as follows
$$R_D(x,y,z)=R_J(x,y,z,z)=\frac32\int_0^\infty\frac{1}{\sqrt{(t+x)(t+y)}(t+z)^\frac32}dt.$$
The linear independence of the symmetric elliptic integrals $R_F(x,y,z)$, $R_D(x,y,z$ and $R_J(x,y,z,p)$, with respect to coefficients that are rational functions, was obtained by Carlson and Gustafson in \cite{CG}.  For the convenience of the reader, we state this as the following lemma.
\vskip 0.2 true cm

\noindent
 {\bf Lemma 2.2} \cite{CG} {\it The functions $R_F(x,y,z)$, $R_D(x,y,z)$, $R_J(x,y,z,p)$, and $\frac{1}{\sqrt{xyz}}$ are
linearly independent with respect to coefficients that are rational functions of $x$, $y$, $z$ and $p$. }\vskip 0.2 true cm

Thanks to Lemma 2.2, we can obtain the linear independence of the complete elliptic integrals $K(k)$, $E(k)$ and $\Pi(\mu,k)$. If $\Pi(\mu,k)$ is regarded as a function of two variables $\mu$ and $k$, where $\mu$ and $k$ are independent of each other, we can obtain the following conclusion.

\vskip 0.2 true cm

\noindent
 {\bf Lemma 2.3} {\it The complete elliptic integrals $K(k)$, $E(k)$ and $\Pi(\mu,k)$ are linearly independent with respect to coefficients that are rational functions of $\mu$ and $k$. }
\vskip 0.2 true cm

\noindent
 {\bf Proof} \ Let $p_1(\mu,k)$, $p_2(\mu,k)$, and $p_3(\mu,k)$  be rational functions of $\mu$ and $k$. We need to prove that
\begin{eqnarray}
p_1(\mu,k)K(k)+p_2(\mu,k)E(k)+p_3(\mu,k)\Pi(\mu,k)\equiv0
\end{eqnarray}
if and only if $p_1(\mu,k)$, $p_2(\mu,k)$ and $p_3(\mu,k)$ are identically $0$.  From formulas 19.25.1 and 19.25.2 in \cite{NIST} or \cite{CG}, one can obtain the relationship between symmetric elliptic integrals and complete elliptic integrals
\begin{eqnarray}\begin{aligned}
&K(k)=R_F(0,1-k^2,1),\\ &K(k)-E(k)=\frac13k^2 R_D(0,1-k^2,1),\\
&\Pi(\mu,k)-K(k)=\frac13 \mu R_J(0,1-k^2,1,1-\mu).
\end{aligned}\end{eqnarray}
Substituting (2.9) into (2.8) gives
\begin{eqnarray}\begin{aligned}
&\big[p_1(\mu,k)+p_2(\mu,k)+p_3(\mu,k)\big]R_F(0,1-k^2,1)\\&
\quad\  -\frac{k^2}{3}p_2(\mu,k)R_D(0,1-k^2,1) +\frac{\mu}{3}p_3(\mu,k)R_J(0,1-k^2,1,1-\mu)\equiv0.
\end{aligned}\end{eqnarray}
It is apparent from Lemma 2.2 and (2.10) that $$p_i(\mu,k)\equiv0,\ i=1,2,3.$$
This completes the proof. $\Box$\vskip 0.2 true cm

In $\Pi(\mu,k)$, provided that $\mu(k)$ is a rational function of $k$, we can still conclude that  $K(k)$, $E(k)$ and $\Pi(\mu,k)$ remain  linearly independent.
\vskip 0.2 true cm

\noindent
 {\bf Lemma 2.4} {\it Let $\mu(k)$ be a rational function of $k$ and $\mu'(0)\neq0$. Then the complete elliptic integrals $K(k)$, $E(k)$ and $\Pi(\mu(k),k)$ are linearly independent with respect to coefficients  that are rational functions of $k$. }
\vskip 0.2 true cm

\noindent
 {\bf Proof}  \   We only need to prove that there exists a point $k_0\in(-1, 1)$ such that the Wronskian determinant $W(k)$ with respect to  $K(k)$, $E(k)$ and $\Pi(\mu(k),k)$ is different from zero at $k_0$, where
 \begin{eqnarray}
 W(k)=\begin{vmatrix}
K(k) & E(k) & \Pi(\mu(k),k) \\
K'(k) & E'(k) & \Pi'(\mu(k),k)\\
K''(k) & E''(k) & \Pi''(\mu(k),k)
\end{vmatrix}.
 \end{eqnarray}
  With the help of Maple, and by substituting (2.2) and (2.3) into (2.11), one can obtain the expression of $W(k)$, which  is quite involved, hence we place it in the Appendix. By using Maple, one has that
 $$\lim\limits_{k\rightarrow0}W(k)=\frac{\pi^3\mu'(0)}{16(1-\mu(0))^\frac32}.$$
Since $\mu'(0)\neq0$ and $\mu(k)\neq1$, one can conclude that there exist a $k_0$ near 0 such that $W(k_0)\neq0.$ This completes the proof. $\Box$\vskip 0.2 true cm

From Lemma 2.4, one can easily obtain the linear independence of the complete elliptic integrals $K(k)$ and $E(k)$.

\vskip 0.2 true cm

\noindent
 {\bf Corollary 2.1} {\it The complete elliptic integrals $K(k)$ and $E(k)$ are linearly independent with respect to coefficients  that are rational functions of $k$. }
 \vskip 0.2 true cm

\noindent
 {\bf Remark 2.1}  If $\Pi(\mu(k),k)$ can be expressed as a linear combination of $K(k)$ and $E(k)$ with rational function coefficients, then
$\mu(k)$ is called a degenerate rational function. As an illustration, $\mu(k)=k^2$ falls into this degenerate category by (2.1), which obviously fails to meet condition $\mu'(0)\neq0$ in Lemma 2.4.
 %\vskip 0.2 true cm

%\noindent
%(ii) It is an open problem whether Lammas 2.2-2.4 is still true if the coefficients are algebraic functions instead of rational functions, see \cite{CG}.
\vskip 0.2 true cm

 \subsection{Proofs of Theorems 1.2-1.4}

With the sufficient preliminary work done, we are now ready to prove the main conclusions of the first part of this paper.
 \vskip 0.2 true cm

\noindent
{\bf Proof of Theorem 1.2}  We begin our proof by eliminating the elliptic integral $\Pi(\mu(k),k)$. To this end, introducing the change of variable
\begin{eqnarray}
Z(k)=X(k)\Pi(\mu(k),k),
\end{eqnarray}
where
\begin{eqnarray*}X(k)=e^{\displaystyle \int\Big(\frac{k}{k^2-\mu}+\frac{(\mu^2-k^2)\mu'(k)}{2\mu(1-\mu)(k^2-\mu)}\Big)dk},\end{eqnarray*}
one can transform (2.2) into
\begin{eqnarray}
\left(  \begin{array}{c}
          K'(k) \\  E'(k)\\ X(k)^{-1}Z'(k)
          \end{array} \right)
=\left(\begin{matrix}
                -\frac{1}{k}&\frac{1}{k(1-k^2)}\\
                 -\frac{1}{k}&\frac{1}{k}\\
                 a_{31}(k)&a_{32}(k)\\
\end{matrix}\right)
\left(
  \begin{array}{c}
          K(k) \\
          E(k)
          \end{array} \right),
 \end{eqnarray}
where $a_{31}(k)$ and $a_{32}(k)$ are given in (2.2).
Denote by $S$ the set  of  zeros of $r(k)$. Then, for any $k\in(-1,1)\setminus S$, an easy calculation gives
\begin{eqnarray}\begin{aligned}
\frac{d}{dk}\Big(\frac{I(k)X(k)}{r(k)}\Big)=&\frac{d}{dk}\Big(\frac{p(k)X(k)}{r(k)}K(k)\Big)+\frac{d}{dk}\Big(\frac{q(k)X(k)}{r(k)}E(k)\Big)+Z'(k)\\
=&\frac{d}{dk}\Big(\frac{p(k)X(k)}{r(k)}\Big)K(k)+\frac{d}{dk}\Big(\frac{q(k)X(k)}{r(k)}\Big)E(k)\\
&+\frac{p(k)X(k)}{r(k)}K'(k)+\frac{q(k)X(k)}{r(k)}E'(k)+Z'(k).
\end{aligned}\end{eqnarray}
It is apparent from (2.13) and (2.14) that
\begin{eqnarray}\begin{aligned}
\frac{d}{dk}\Big(\frac{I(k)X(k)}{r(k)}\Big)=&r^{-2}(k)\Big[p'(k)r(k)X(k)+p(k)r(k)X'(k)-p(k)r'(k)X(k)\\
& -\frac{p(k)r(k)X(k)}{k}-\frac{q(k)r(k)X(k)}{k}+a_{31}(k)r^2(k)X(k)\Big]K(k)\\
&+r^{-2}(k)\Big[q'(k)r(k)X(k)+q(k)r(k)X'(k)-q(k)r'(k)X(k)\\
& +\frac{p(k)r(k)X(k)}{k(1-k^2)}+\frac{q(k)r(k)X(k)}{k}+a_{32}(k)r^2(k)X(k)\Big]E(k)\\
:=&r^{-2}(k)\big[M(k)K(k)+N(k)E(k)\big],
\end{aligned}\end{eqnarray}
where $M(k)$ and $N(k)$ are functions of $k$. This completes the proof. $\Box$ \vskip 0.2 true cm

In the sequel, we will use the notation $\#\{k\in\Sigma|\varphi(k)=0\}$ to indicate the number of isolated zeros of the real function $\varphi(k)$ in the set $\Sigma$, taking into account their multiplicities.

 \vskip 0.2 true cm

\noindent
{\bf Proof of Theorem 1.3}\ \, Suppose that $p(k)$, $q(k)$, $r(k)$ and  $\mu(k)$ are real polynomials of degree $m$, $n$, $l$ and $s$, respectively. On account of Theorem 1.2, it suffices to estimate the number of zeros of $M(k)K(k)+N(k)E(k)$. Since $p(k)$, $q(k)$, $r(k)$ and  $\mu(k)$ are polynomials, this enables us to accomplish this.

When  $s\geq2$ and $k\in(-1,1)$, by substituting $X(k)$ in (2.12) and $p(k)$, $q(k)$, $r(k)$ and  $\mu(k)$ into (2.15),
%\begin{eqnarray}\begin{aligned}
%\frac{d}{dk}\Big(\frac{I(k)X(k)}{r(k)}\Big)%=&\frac{d}{dk}\Big(\frac{p(k)X(k)}{r(k)}K(k)\Big)+\frac{d}{dk}\Big(\frac{q(k)X(k)}{r(k)}E(k)\Big)+\frac{dZ(k)}{dk}\\
%=&\frac{1}{r^2(k)}\Big[\big(p'(k)X(k)r(k)+ p(k)X'(k)r(k)-p(k)X(k)r'(k)\big)K(k)\\&+\big(q'(k)X(k)r(k)+ q(k)X'(k)r(k)-q(k)X(k)r'(k)\big)E(k)\\&+p(k)X(k)r(k)K'(k)+p(k)X(k)r(k)E'(k)+Z'(k)\Big].
%\end{aligned}\end{eqnarray}
 one finds that
\begin{eqnarray}\begin{aligned}
M(k)K(k)+N(k)E(k)=&\frac{X(k)\big[M_1(k)K(k)+N_1(k)E(k)\big]}{k(1-k^2)\mu(k)(1-\mu(k))\big(k^2-\mu(k)\big)}\\
:=&\frac{X(k)Q_1(k)}{k(1-k^2)r^2(k)\mu(k)(1-\mu(k))\big(k^2-\mu(k)\big)},
\end{aligned}\end{eqnarray}
where %$M_1(k)$ and $N_1(k)$ are polynomials of $k$ and
\begin{eqnarray}\begin{aligned}
M_1(k)=&(k-k^3)\big(k^2-\mu(k)\big)\big(\mu(k)-\mu^2(k)\big)[p'(k)r(k)-p(k)r'(k)]\\
&+(k-k^3)\big[k\big(\mu(k)-\mu^2(k)\big)+\frac12\big(\mu^2(k)-k^2\big)\mu'(k)\big]p(k)r(k)\\
&-(1-k^2)\big(k^2-\mu(k)\big)\big(\mu(k)-\mu^2(k)\big)\big[p(k)+q(k)\big]r(k)\\
&-\frac12(k-k^3)\big(k^2-\mu(k)\big)\mu'(k)r^2(k),\\
N_1(k)=&(k-k^3)\big(k^2-\mu(k)\big)\big(\mu(k)-\mu^2(k)\big)[q'(k)r(k)-q(k)r'(k)]\\
&+(k-k^3)\big[k\big(\mu(k)-\mu^2(k)\big)+\frac12\big(\mu^2(k)-k^2\big)\mu'(k)\big]q(k)r(k)\\
&+\big(k^2-\mu(k)\big)\big(\mu(k)-\mu^2(k)\big)\big[p(k)+(1-k^2)q(k)\big]r(k)\\
&+k\mu(k)\big[k\big(1-\mu(k)\big)-\frac12(1-k^2)\mu'(k)\big]r^2(k),\\
\end{aligned}\end{eqnarray}
with
\begin{eqnarray}\begin{aligned}
&\deg M_1(k)\leq\max\{n+l+3s+2,\,m+l+3s+2,\,2l+2s+2\},\\ &\deg N_1(k)\leq\max\{n+l+3s+2,\,m+l+3s,\,2l+2s+2\}.
\end{aligned}\end{eqnarray}
Using the Theorem 1 in \cite{GLLZ} (for the reader's convenience, we restate this theorem as Theorem 1.1 in this paper), one derives that
\begin{eqnarray}
\begin{aligned}
\#\{k\in(-1,1)| Q_1(k)=0\}\leq \deg M_1(k)+\deg N_1(k)+2\\
\leq\begin{cases}
2\max\{m,n\}+2l+6s+6, \ \textup{if}\ l\leq\max\{m,n\}+s,\\
4l+4s+6, \qquad\qquad\qquad\ \ \ \textup{if}\ l>\max\{m,n\}+s,
\end{cases}
\end{aligned}
\end{eqnarray}
in view of (2.18).
Since the function $$\frac{X(k)}{k(1-k^2)\mu(k)(1-\mu(k))\big(k^2-\mu(k)\big)}$$ is non-vanishing  on the interval $(-1,1)$, it is apparent from (2.15), (2.16) and (2.19) that
\begin{eqnarray*}
\begin{aligned}
&\#\{k\in(-1,1)| I(k)=0\}\leq\deg Q_1(k)+\deg r(k)+1\\
&\qquad\leq \begin{cases}
2\max\{m,n\}+3l+6s+7,\ \textup{if}\ l\leq \max\{m,n\}+s,\\
5l+4s+7,\qquad\qquad\qquad\ \ \, \, \textup{if}\ l>\max\{m,n\}+s.
\end{cases}
\end{aligned}
\end{eqnarray*}

When $s=1$, similar to (2.16), one concludes that
\begin{eqnarray*}\begin{aligned}
M(k)K(k)+N(k)E(k)=&\frac{X(k)\big[M_1(k)K(k)+N_1(k)E(k)\big]}{k(1-k^2)\mu(k)(1-\mu(k))\big(k^2-\mu(k)\big)},
%:=&\frac{X(k)W_1(k)}{k(1-k^2)r^2(k)\mu(k)(1-\mu(k))\big(k^2-\mu(k)\big)},
\end{aligned}\end{eqnarray*}
where $M_1(k)$ and $N_1(k)$ are defined as in (2.17) and
\begin{eqnarray*}\begin{aligned}
&\deg M_1(k)\leq\max\{n+l+6,\,m+l+6,\,2l+5\},\\ &\deg N_1(k)\leq\max\{n+l+6,\,m+l+4,\,2l+4\}.
\end{aligned}\end{eqnarray*}
It follows that
\begin{eqnarray*}
\begin{aligned}
\#\{k\in(-1,1)| I(k)=0\}
\leq \begin{cases}
2\max\{m,n\}+3l+15,\ \textup{if}\ l\leq \max\{m,n\}+1,\\
5l+13,\qquad\qquad\qquad\ \ \, \, \textup{if}\ l>\max\{m,n\}+1.
\end{cases}
\end{aligned}
\end{eqnarray*}

When $s=0$, that is, $\mu(k)$ is a nonzero constant independent of $k$. It follows that
$$\mu'(k)=0,\ Z(k)=\sqrt{|k^2-\mu|}\Pi(\mu(k),k).$$
A similar calculation leads to
\begin{eqnarray*}\begin{aligned}
M(k)K(k)+N(k)E(k)=&\frac{M_1(k)K(k)+N_1(k)E(k)}{k(1-k^2)\sqrt{|k^2-\mu|}},
%:=&\frac{X(k)W_1(k)}{k(1-k^2)r^2(k)\mu(k)(1-\mu(k))\big(k^2-\mu(k)\big)},
\end{aligned}\end{eqnarray*}
where
\begin{eqnarray*}\begin{aligned}
M_1(k)=&(k-k^3)\big(\mu-k^2\big)[p'(k)r(k)-p(k)r'(k)]-\mu(1-k^2) p(k)r(k)\\
&-(1-k^2)\big(\mu-k^2\big)q(k)r(k),\\
N_1(k)=&(k-k^3)\big(\mu-k^2\big)[q'(k)r(k)-q(k)r'(k)]-\mu\big(\mu-k^2\big) p(k)r(k)\\
&+(1-k^2)\big(\mu-2k^2\big)q(k)r(k)-k^2r^2(k),\\
\end{aligned}\end{eqnarray*}
with
\begin{eqnarray*}\begin{aligned}
&\deg M_1(k)\leq\max\{n+l+4,\,m+l+4\},\\ &\deg N_1(k)\leq\max\{n+l+4,\,m+l+2,\,2l+2\}.
\end{aligned}\end{eqnarray*}
Therefore, on can obtain that
\begin{eqnarray*}
\begin{aligned}
\#\{k\in(-1,1)| I(k)=0\}
\leq \begin{cases}
2\max\{m,n\}+3l+11,\ \textup{if}\ l\leq \max\{m,n\}+2,\\
\max\{m,n\}+4l+9,\ \ \ \, \ \textup{if}\ l>\max\{m,n\}+2.
\end{cases}
\end{aligned}
\end{eqnarray*}
%The remainder of the argument is analogous to that in the case $s\geq1$ and is considerably simpler, hence we omit the proof.

We are now in a position to prove the second statement in this theorem. By Lemma 2.4, one gets that there exist non-zero  polynomials $p(k)$, $q(k)$ and $r(k)$ such that $I(k)\not\equiv0$. For any given $k_i\in(-1,1)$, $i=1,2,\cdots, m+n+l+2$, consider the following system of linear equations
\begin{eqnarray}\begin{aligned}
\sum\limits_{j=0}^ma_jk^j_iK(k_i)+\sum\limits_{j=0}^nb_jk^j_iE(k_i)+\sum\limits_{j=0}^lc_jk^j_i\Pi\big(\mu(k_i),k_i\big)=0,\\ i=1,2,\cdots,m+n+l+2.
\end{aligned}\end{eqnarray}
Since  system (2.20) contains $m+n+l+3$ unknown variables $a_i,i=0,1,\cdots m$,   $b_i,i=0,1,\cdots n$ and $c_i,i=0,1,\cdots l$, and only $m+n+l+2$ equations, this implies that it has at least one nontrivial solution
 $$\{a_0,a_1,\cdots a_m, b_0,b_1,\cdots,b_n,c_0,c_1,\cdots,c_l\}.$$
This completes the proof. $\Box$

 \vskip 0.2 true cm

\noindent
{\bf Proof of Theorem 1.4}\ \, When $\mu(k)=\frac{2k^2}{1+k^2}$, the matrix $\mathbf{A}(k)$ in the Picard-Fuchs (2.2) can be written as
\begin{eqnarray*}
\mathbf{A}(k)=\left(\begin{matrix}
                -\frac{1}{k}&\frac{1}{k(1-k^2)}&0\\
                 -\frac{1}{k}&\frac{1}{k}&0\\
                 \frac{1}{k(k^2-1)}&\frac{1}{k(1-k^2)}&\frac{k(k^2+3)}{1-k^4}\\
\end{matrix}\right).
\end{eqnarray*}
Hence,  by the substitution
\begin{eqnarray}
Z(k)=X(k)\Pi(\mu(k),k)=\frac{1-k^2}{\sqrt{1+k^2}}\Pi(\mu(k),k),
\end{eqnarray}
one finds that
\begin{eqnarray}
k(1-k^2)\left(  \begin{array}{c}
          K'(k) \\  E'(k)\\ \frac{\sqrt{1+k^2}}{1-k^2}Z'(k)
          \end{array} \right)
=\left(\begin{matrix}
                k^2-1&1\\
                 k^2-1&1-k^2\\
                 -1&1\\
\end{matrix}\right)
\left(
  \begin{array}{c}
          K(k) \\
          E(k)
          \end{array} \right).
 \end{eqnarray}
 It follows from (2.21) and (2.22) that
\begin{eqnarray}\begin{aligned}
\frac{d}{dk}\Big(\frac{1-k^2}{r(k)\sqrt{1+k^2}}I(k)\Big)=&\frac{M_2(k)K(k)+N_2(k)E(k)}{r^2(k)k(1+k^2)^\frac32},\ k\in(-1,1)\setminus S,
\end{aligned}\end{eqnarray}
where
$$\begin{aligned}
M_2(k)=&(k-k^5)[p'(k)r(k)-p(k)r'(k)]-(k^4+3k^2)p(k)r(k)\\&-(1-k^4)p(k)r(k)-(1-k^4)q(k)r(k)-(1+k^2)r^2(k),\\
N_2(k)=&(k-k^5)[q'(k)r(k)-q(k)r'(k)]-(k^4+3k^2)q(k)r(k)\\&+(1+k^2)p(k)r(k)+(1-k^4)q(k)r(k)+(1+k^2)r^2(k)
\end{aligned}$$
with%$M_2(k)$ and $N_2(k)$ are polynomials of $k$ and
\begin{eqnarray}\begin{aligned}
&\deg M_2(k)\leq\max\{n+l+4,\,m+l+4,\,2l+2\},\\ &\deg N_2(k)\leq\max\{n+l+4,\,m+l+2,\,2l+2\}.
\end{aligned}\end{eqnarray}
Therefore, using the Theorem 1 in \cite{GLLZ}, (2.23) and (2.24), one has that
\begin{eqnarray*}
\begin{aligned}
\#\{k\in(-1,1)| I(k)=0\}\leq&\deg M_2(k)+\deg N_2(k)+\deg r(k)+3\\
\leq &\begin{cases}
2\max\{m,n\}+3l+11,\ \textup{if}\ l\leq \max\{m,n\}+2,\\
5l+7,\qquad\qquad\qquad\quad \,\ \textup{if}\ l>\max\{m,n\}+2,\\
\end{cases}
\end{aligned}
\end{eqnarray*}
together with Rolle's Theorem. The proof of the second part is identical to that of Theorem 1.3. This completes the proof. $\Box$

\section{Perturbation of a Hamiltonian triangle}
 \setcounter{equation}{0}
\renewcommand\theequation{3.\arabic{equation}}

 When $\varepsilon=0$, the corresponding Hamiltonian function of system (1.6) is
\begin{eqnarray*}
H(x,y)=x^2y(1-x-y).
\end{eqnarray*}
Furthermore, system $(1.6)|_{\varepsilon=0}$ has three invariant straight lines $x=0$, $y=0$ and $y=-x+1$.  An elementary center $(\frac{1}{2},\frac{1}{4})$ of it lies inside the triangle formed by these lines, while a hyperbolic saddle point $(1, 0)$ is one vertex of the triangle. System $(1.6)|_{\varepsilon=0}$ has a family of periodic orbits surrounded the elementary center $(\frac{1}{2},\frac{1}{4})$ defined by
$$\Gamma_h=\{(x,y)|H(x,y)=h,h\in(0,\frac{1}{64})\}.$$
%It is easy to see that  $x=0$ is also a singular straight line.
The straight line $y=-\frac{1}{2}x+\frac12$ divides the periodic orbit $\Gamma_h$ into upper and lower arcs, which are explicitly given by
$$\Gamma^+_h=\{(x,y)|y=-\frac{1}{2}x+\frac{1}{2}+\frac{1}{2}x^{-1}\sqrt{x^4-2x^3+x^2-4h}\},$$ and $$\Gamma^-_h=\{(x,y)|y=-\frac{1}{2}x+\frac{1}{2}-\frac{1}{2}x^{-1}\sqrt{x^4-2x^3+x^2-4h}\}.$$
The coordinates of intersection points between this straight line and the closed orbit $\Gamma_h$ are
$$ A\Big(\frac12-\frac12\sqrt{1-8\sqrt{h}},\ \frac14+\frac14\sqrt{1-8\sqrt{h}}\Big),$$and $$B\Big(\frac12+\frac12\sqrt{1-8\sqrt{h}},\ \frac14-\frac14\sqrt{1-8\sqrt{h}}\Big).$$
When $h$ approaches to $0^+$, $\Gamma_h$ approaches to  the triangle with vertices $(0,0),$ $(1,0)$ and $(0,1)$. When $h$ approaches to $\frac{1}{64}^-$, $\Gamma_h$ approaches to the elementary center $(\frac{1}{2},\frac{1}{4})$,  as shown in Fig. 1.

 \subsection{Algebraic structure of the Melnikov function}
 \setcounter{equation}{0}

In order to estimate the number of zeros of the first order Melnikov function $I(h)$, one should study the algebraic structure of $I(h)$. To this end, we denote
\begin{eqnarray}\begin{aligned}
&I_{i,j}(h)=\frac{1}{i+1}\int_{\Gamma_h^+}x^{i+1}y^jdy,\\ &J_{i,j}(h)=\frac{1}{i+1}\int_{\Gamma_h^-}x^{i+1}y^jdy,\ h\in(0,\frac{1}{64}).
\end{aligned}\end{eqnarray}
In the rest of this paper, we always denote by  $\phi_n(\cdot)$ an odd real polynomial of degree no more than $n$ and by  $P_m(\cdot)$ a real polynomial of degree no more than $m$ unless the opposite is claimed. Note that this is just a symbol. For example, $\phi_{n-2}(\cdot)+\phi_n(\cdot)=\phi_n(\cdot)$ and $P_{m-1}(\cdot)+P_m(\cdot)=P_m(\cdot)$. We first prove that the first order Melnikov function $I(h)$ admits a representation as a linear combination of generator integrals with polynomial coefficients. The following lemma helps us attain our goal.
\vskip 0.2 true cm

\noindent
 {\bf Lemma 3.1} {\it The following relationships hold:
\begin{eqnarray}\begin{aligned}
\int_{\Gamma^+_{h}}x^iy^jdx=\begin{cases}-\frac{j}{i+1}\int_{\Gamma^+_{h}}x^{i+1}y^{j-1}dy+\phi^+_{i+j+\frac{1+(-1)^{i+j}}{2}}(\sqrt{1-8\sqrt{h}}),\ j\neq i+1,\\
-\frac{j}{i+1}\int_{\Gamma^+_{h}}x^{i+1}y^{j-1}dy,\qquad\qquad\qquad\qquad\qquad\qquad\, \, j=i+1,\end{cases}\\
\int_{\Gamma^-_{h}}x^iy^jdx=
\begin{cases}
-\frac{j}{i+1}\int_{\Gamma^-_{h}}x^{i+1}y^{j-1}dy+\phi^-_{i+j+\frac{1+(-1)^{i+j}}{2}}(\sqrt{1-8\sqrt{h}}),\ j\neq i+1,\\
-\frac{j}{i+1}\int_{\Gamma^-_{h}}x^{i+1}y^{j-1}dy,\qquad\qquad\qquad\qquad\qquad\qquad\, \, j=i+1,\end{cases}
\end{aligned}\end{eqnarray}
where $\phi^\pm_{l}(\cdot)$ are odd polynomials of degree no more than $l$.}
\vskip 0.2 true cm

\noindent
{\bf Proof}\, By a straightforward computation using the Green's Formula, one derives that
\begin{eqnarray}\begin{aligned}
\int_{\Gamma^+_{h}}x^iy^jdx=&\int_{\Gamma^+_{h}\cup \overrightarrow{AB}}x^{i}y^{j}dy-\int_{\overrightarrow{AB}}x^{i}y^{j}dy\\
=&-j\iint\limits_\Omega x^iy^{j-1}dxdy-2^{-j}\int_{x_0}^{x_1}x^{i}(1-x)^{j}dx,\\
\end{aligned}\end{eqnarray}
where $\Omega$ is the interior region bounded by $\Gamma^+_{h}$ and $\overrightarrow{AB}$, and
$$x_0=\frac12-\frac12\sqrt{1-8\sqrt{h}},\ x_1=\frac12+\frac12\sqrt{1-8\sqrt{h}}.$$
Similarly, one can obtain
\begin{eqnarray}\begin{aligned}
\int_{\Gamma^+_{h}}x^iy^jdy=i\iint\limits_\Omega x^{i-1}y^{j}dxdy+2^{-j-1}\int_{x_0}^{x_1}x^{i}(1-x)^{j}dx.
\end{aligned}\end{eqnarray}
An easy maniputation using (3.3) and (3.4) gives
\begin{eqnarray}\begin{aligned}
\int_{\Gamma^+_{h}}x^iy^jdx=&-\frac{j}{i+1}\int_{\Gamma^+_{h}}x^{i+1}y^{j-1}dy-\frac{1}{2^j}\int_{x_0}^{x_1}x^{i}(1-x)^{j}dx\\&+\frac{j}{(i+1)2^j}\int_{x_0}^{x_1}x^{i+1}(1-x)^{j-1}dx.
\end{aligned}\end{eqnarray}
To facilitate subsequent discussion, we set $$F(u)=\frac{j}{(i+1)}\int_{x_0}^{x_1}x^{i+1}(1-x)^{j-1}dx-\int_{x_0}^{x_1}x^{i}(1-x)^{j}dx,$$
where \begin{eqnarray}u=\sqrt{1-8\sqrt{h}},\ u\in(0,1).\end{eqnarray}
When $j=i+1$, the change of variables $t=x-\frac12$ allows us to obtain
\begin{eqnarray}\begin{aligned}
F(u)=&\int_{x_0}^{x_1}x^{i}(1-x)^{i}(2x-1)dx\\
=&2\int_{-\frac u2}^{\frac u2}t(\frac14-t^2)^{i}dt.
\end{aligned}\end{eqnarray}
Due to the integrand being an odd function in the integral (3.7), it equals zero immediately by symmetry. When $j\neq i+1$,  by the same transformation $t=x-\frac12$, one has that
$$\begin{aligned}
F(u)=\frac{j}{i+1}\int_{-\frac u2}^{\frac u2}(t+\frac12)^{i+1}(\frac12-t)^{j-1}dt-\int_{-\frac u2}^{\frac u2}(t+\frac12)^{i}(\frac12-t)^{j}dt.
\end{aligned}$$
It is clear that $F(u)+F(-u)=0$. Therefore, $F(u)$ is an odd polynomial of $u$, and its highest degree is $i+j+\frac{1+(-1)^{i+j}}{2}$. This indicates that the first relationship in (3.2) is valid.
The second one can be proved similarly. The proof is complete. $\Box$
\vskip 0.2 true cm

\noindent
 {\bf Remark 3.1}  In what follows, we will not distinguish between the cases $j\neq i+1$ and $j=i+1$, as this distinction does not affect the expression of the first order Melnikov function $I(h)$. This is because the structure of $I(h)$ is entirely determined by the case $j\neq i+1$.
\vskip 0.2 true cm

\noindent
 {\bf Lemma 3.2} {\it For $h\in(0,\frac{1}{64})$ and any natural number $n\geq3$, it holds that
 \begin{eqnarray}\begin{aligned}
I(h)=&\Big(\sum\limits_{i=0}^{[\frac{n}{4}]}\alpha^+_ih^i\Big)I_{0,0}(h)+\Big(\sum\limits_{i=0}^{[\frac{n-1}{4}]}\beta^+_ih^i\Big)I_{2,0}(h)
+\Big(\sum\limits_{i=0}^{[\frac{n-3}{4}]}\gamma^+_ih^i\Big)I_{3,0}(h)\\
&+\Big(\sum\limits_{i=0}^{[\frac{n}{4}]}\alpha^-_ih^i\Big)J_{0,0}(h)+\Big(\sum\limits_{i=0}^{[\frac{n-1}{4}]}\beta^-_ih^i\Big)J_{2,0}(h)
+\Big(\sum\limits_{i=0}^{[\frac{n-3}{4}]}\gamma^-_ih^i\Big)J_{3,0}(h)\\
&++h\phi_{n-3+\frac{1-(-1)^n}{2}}(\sqrt{1-8\sqrt{h}})+\phi_{n+1+\frac{1-(-1)^n}{2}}(\sqrt{1-8\sqrt{h}}),
\end{aligned}\end{eqnarray}
where $\alpha^\pm_i$, $\beta^\pm_i$ and $\gamma^\pm_i$ are constants and $\phi_{k}(\cdot)$ is an odd polynomial of degree no more than $k$.}
\vskip 0.2 true cm

\noindent
{\bf Proof.}\, It follows from (3.2) that
\begin{eqnarray}
\begin{aligned}
I(h)=&\sum\limits_{i+j=0}^{n}b^+_{i,j}\int_{\Gamma^+_h}x^iy^jdx-\sum\limits_{i+j=0}^{n}a^+_{i,j}\int_{\Gamma^+_h}x^iy^jdy\\
&+\sum\limits_{i+j=0}^{n}b^-_{i,j}\int_{\Gamma^-_h}x^iy^jdx-\sum\limits_{i+j=0}^{n}a^-_{i,j}\int_{\Gamma^-_h}x^iy^jdy\\
\triangleq&\sum\limits_{i+j=0}^{n}\xi_{i,j}I_{i,j}(h)+\sum\limits_{i+j=0}^{n}\eta_{i,j}J_{i,j}(h)+\phi_{n+\frac{1+(-1)^n}{2}}(\sqrt{1-8\sqrt{h}}),
\end{aligned}
\end{eqnarray}
where $\xi_{i,j}$ and $\eta_{i,j}$ are constants.

Now we claim that
\begin{eqnarray}\begin{aligned}
\sum\limits_{i+j=0}^{n}\xi_{i,j}I_{i,j}(h)=&\Big(\sum\limits_{i=0}^{[\frac{n}{4}]}\alpha^+_ih^i\Big) I_{0,0}(h)+\Big(\sum\limits_{i=0}^{[\frac{n-1}{4}]}\beta^+_ih^i\Big)I_{2,0}(h)+\Big(\sum\limits_{i=0}^{[\frac{n-3}{4}]}\gamma^+_ih^i\Big)I_{3,0}(h)
\\&+h\phi^+_{n-3+\frac{1-(-1)^n}{2}}(\sqrt{1-8\sqrt{h}})+\phi^+_{n+1+\frac{1-(-1)^n}{2}}(\sqrt{1-8\sqrt{h}}),
\end{aligned}\end{eqnarray}
where $\alpha^+_i$, $\beta^+_i$ and $\gamma^+_i$ are constants, and $\phi^+_{l}(\cdot)$ is an odd polynomials of degree no more than $l$.

In fact, an easy calculation using $H(x,y)=h$ yields
\begin{eqnarray}\begin{aligned}
&h(i+1)I_{i,j}(h)=(i+3)I_{i+2,j+1}(h)-(i+4)I_{i+3,j+1}(h)-(i+3)I_{i+2,j+2}(h).
\end{aligned}\end{eqnarray}
It follows from (3.2) that
\begin{eqnarray}
I_{i,j}(h)=%\begin{cases}
-\frac{1}{j+1}\int_{\Gamma^+_h}x^iy^{j+1}dx+\frac{1}{j+1}\phi^+_{i+j+1+\frac{1-(-1)^{i+j}}{2}}(\sqrt{1-8\sqrt{h}}).%\ i\neq j,\\
%-\frac{1}{j+1}\int_{\Gamma^+_h}x^iy^{j+1}dx,\qquad\qquad\qquad\qquad\qquad\qquad\qquad\ \ \ i=j.
%\end{cases}
\end{eqnarray}
Multiplying both sides of (3.12) by $h(j+1)$ implies
\begin{eqnarray}\begin{aligned}
h(j+1)I_{i,j}(h)=&(j+2)I_{i+2,j+1}(h)-(j+2)I_{i+3,j+1}(h)-(j+3)I_{i+2,j+2}(h)\\&+h\phi^+_{i+j+1+\frac{1-(-1)^{i+j}}{2}}(\sqrt{1-8\sqrt{h}})+\phi^+_{i+j+5+\frac{1-(-1)^{i+j}}{2}}(\sqrt{1-8\sqrt{h}}).
\end{aligned}\end{eqnarray}
Therefore, by multiplying (3.11) by $j+3$ and adding (3.13) multiplied by $-(i+3)$, one can obtain that
\begin{eqnarray}\begin{aligned}
I_{i+3,j+1}(h)=&\frac{1}{i+j+6}\big[2h(j-i)I_{i,j}(h)+(i+3)I_{i+2,j+1}(h)\big]
\\&-\frac{i+3}{i+j+6}\Big[h\phi^+_{i+j+1+\frac{1-(-1)^{i+j}}{2}}(\sqrt{1-8\sqrt{h}})\\&+\phi^+_{i+j+5+\frac{1-(-1)^{i+j}}{2}}(\sqrt{1-8\sqrt{h}}\Big].
\end{aligned}\end{eqnarray}
Similarly, one concludes that
\begin{eqnarray}\begin{aligned}
I_{i+2,j+2}(h)=&\frac{1}{i+j+6}\big[h(i-3j-2)I_{i,j}(h)+(j+2)I_{i+2,j+1}(h)\big]
\\&+\frac{i+4}{i+j+6}\Big[h\phi^+_{i+j+1+\frac{1-(-1)^{i+j}}{2}}(\sqrt{1-8\sqrt{h}})\\&+\phi^+_{i+j+5+\frac{1-(-1)^{i+j}}{2}}(\sqrt{1-8\sqrt{h}}\Big],
\end{aligned}\end{eqnarray}
and
\begin{eqnarray}
\begin{aligned}
I_{i+2,j+1}(h)=&\frac{1}{i-2j-1}\big[(i-3j-2)I_{i+3,j+1}(h)+2(i-j)I_{i+2,j+2}(h)\big]
\\&-\frac{i+1}{i-2j-1}\Big[h\phi^+_{i+j+1+\frac{1-(-1)^{i+j}}{2}}(\sqrt{1-8\sqrt{h}})\\&+\phi^+_{i+j+5+\frac{1-(-1)^{i+j}}{2}}(\sqrt{1-8\sqrt{h}}\Big].
\end{aligned}
\end{eqnarray}

We are now in a position to prove the claim by induction on $n$. By taking $(i, j) = (-1, -1)$ in (3.14) and (3.15) respectively, one obtains that
\begin{eqnarray}\begin{aligned}
&I_{1,0}(h)=2I_{2,0}(h)+\phi^+_{3}(\sqrt{1-8\sqrt{h}}),\\
&I_{1,1}(h)=\frac{1}{2}I_{2,0}(h)+\phi^+_{3}(\sqrt{1-8\sqrt{h}}).
\end{aligned}
\end{eqnarray}
By substituting $(i,j) = (-2,-1), (-2,0),(0,-1),(-1,0)$ and $(-2,1)$ into (3.16) respectively, one concludes that
\begin{eqnarray}\begin{aligned}
&I_{0,1}(h)=\frac{1}{2}I_{0,0}(h)-\frac{1}{2}I_{1,0}(h)+\frac12\phi^+_{3}(\sqrt{1-8\sqrt{h}}),\\
&I_{0,2}(h)=\frac{3}{4}I_{0,1}(h)-I_{1,1}(h)+\frac14\phi^+_{3}(\sqrt{1-8\sqrt{h}}),\\
&I_{2,1}(h)=\frac{1}{2}I_{2,0}(h)-\frac{1}{2}I_{3,0}(h)+\frac12\phi^+_{5}(\sqrt{1-8\sqrt{h})}+\frac12 h\phi^+_{1}(\sqrt{1-8\sqrt{h}}),\\
&I_{1,2}(h)=I_{1,1}(h)-\frac{3}{2}I_{2,1}(h),\\
&I_{0,3}(h)=\frac{5}{6}I_{0,2}(h)-\frac{7}{6}I_{1,2}(h)+\frac16\phi^+_{5}(\sqrt{1-8\sqrt{h})}+\frac16 h\phi^+_{1}(\sqrt{1-8\sqrt{h}}).
\end{aligned}\end{eqnarray}
Hence, for $n=3$, a direct calculation using (3.17) and (3.18) yields
$$\begin{aligned}\sum\limits_{i+j=0}^3\xi_{i,j}I_{i,j}(h)=&\alpha^+_0I_{0,0}(h)+\beta^+_0I_{2,0}(h)+\gamma^+_0I_{3,0}(h)\\&+h\phi^+_{1}(\sqrt{1-8\sqrt{h}})+\phi^+_{5}(\sqrt{1-8\sqrt{h}}),\end{aligned}$$
where $\alpha^+_0,\beta^+_0$ and $\gamma^+_0$ are constants. That is, the claim is valid for $n=3$.

Now assume that the claim holds for all $i+j\leq n$. Then, by taking $(i,j)=(0,n-3)$ in (3.15) and $(i,j)=(0,n-3),(1,n-4),\cdots,(n-3,0)$ in (3.14), respectively, one has that
\begin{eqnarray}
\begin{cases}
I_{2,n-1}(h)=&\frac{1}{n+3}[(3n-7)hI_{0,n-3}(h)+(n-1)I_{2,n-2}(h)]\\&+\frac{4}{n+3}\big[h\phi^+_{n-2+\frac{1+(-1)^n}{2}}(\sqrt{1-8\sqrt{h}})+\phi^+_{n+2+\frac{1+(-1)^n}{2}}(\sqrt{1-8\sqrt{h}})\big],\\
I_{3,n-2}(h)=&\frac{1}{n+3}[2(n-3)hI_{0,n-3}(h)+3I_{2,n-2}(h)]\\&-\frac{3}{n+3}\big[h\phi^+_{n-2+\frac{1+(-1)^n}{2}}(\sqrt{1-8\sqrt{h}})+\phi^+_{n+2+\frac{1+(-1)^n}{2}}(\sqrt{1-8\sqrt{h}})\big],\\
I_{4,n-3}(h)=&\frac{1}{n+3}[2(n-5)hI_{1,n-4}(h)+4I_{3,n-3}(h)]\\&-\frac{4}{n+3}\big[h\phi^+_{n-2+\frac{1+(-1)^n}{2}}(\sqrt{1-8\sqrt{h}})+\phi^+_{n+2+\frac{1+(-1)^n}{2}}(\sqrt{1-8\sqrt{h}})\big],\\
   \qquad\qquad\              \vdots\\
I_{n-1,2}(h)=&\frac{1}{n+3}[2(5-n)hI_{n-4,1}(h)+(n-1)I_{n-2,2}(h)]\\&-\frac{n-1}{n+3}\big[h\phi^+_{n-2+\frac{1+(-1)^n}{2}}(\sqrt{1-8\sqrt{h}})+\phi^+_{n+2+\frac{1+(-1)^n}{2}}(\sqrt{1-8\sqrt{h}})\big],\\
I_{n,1}(h)=&\frac{1}{n+3}[2(3-n)hI_{n-3,0}(h)+nI_{n-1,1}(h)]\\&-\frac{n}{n+3}\big[h\phi^+_{n-2+\frac{1+(-1)^n}{2}}(\sqrt{1-8\sqrt{h}})+\phi^+_{n+2+\frac{1+(-1)^n}{2}}(\sqrt{1-8\sqrt{h}})\big].
\end{cases}\end{eqnarray}
By substituting $(i,j)=(-2,n-1),(-1,n-2)$ and $(n-2,-1)$ into (3.16), one obtains that
\begin{eqnarray}
\begin{cases}
I_{0,n+1}(h)=&\frac{1}{2n+2}[(2n+1)I_{0,n}(h)-(3n+1)I_{1,n}(h)]\\&+\frac{1}{2n+2}\big[h\phi^+_{n-2+\frac{1+(-1)^n}{2}}(\sqrt{1-8\sqrt{h}})+\phi^+_{n+2+\frac{1+(-1)^n}{2}}(\sqrt{1-8\sqrt{h}})\big],\\
I_{1,n}(h)=&I_{1,n-1}(h)-\frac{3}{2}I_{2,n-1}(h),\\
I_{n+1,0}(h)=&I_{n,0}(h)-2I_{n,1}(h)\\&+h\phi^+_{n-2+\frac{1+(-1)^n}{2}}(\sqrt{1-8\sqrt{h}})+\phi^+_{n+2+\frac{1+(-1)^n}{2}}(\sqrt{1-8\sqrt{h}}).
\end{cases}\end{eqnarray}
In order to utilize the induction hypothesis, it is necessary to remove the terms with subscript sums equal to $n+1$  in the right hand side of equation (3.20) by employing the first and last relations provided in (3.19). A direct calculation shows that (3.20) can be transformed into the following form
\begin{eqnarray}
\begin{cases}
I_{0,n+1}(h)=&\frac{1}{4(n+1)(n+3)}\big[(9n^2-6n-3)I_{2,n-2}(h)-(6n^2+20n+6)I_{1,n-1}(h)\\&-(27n^2-54n-21)hI_{0,n-3}(h)+(n+3)(4n+2)I_{0,n}(h)\big]
\\&+\frac{19n+9}{2(n+1)(n+3)}\big[h\phi^+_{n-2+\frac{1+(-1)^n}{2}}(\sqrt{1-8\sqrt{h}})+\phi^+_{n+2+\frac{1+(-1)^n}{2}}(\sqrt{1-8\sqrt{h}})\big],\\
I_{1,n}(h)=&\frac{1}{2n+6}\big[(2n+6)I_{1,n-1}(h)-(3n-3)I_{2,n-2}(h)+(9n-27)hI_{0,n-3}(h)\big]\\&-\frac{6}{n+3}\big[h\phi^+_{n-2+\frac{1+(-1)^n}{2}}(\sqrt{1-8\sqrt{h}})+\phi^+_{n+2+\frac{1+(-1)^n}{2}}(\sqrt{1-8\sqrt{h}})\big],\\
I_{n+1,0}(h)=&\frac{1}{n+3}\big[(n+3)I_{n,0}(h)+4(n-3)hI_{n-3,0}(h)-2nI_{n-1,1}(h)\big]\\&+\frac{3n+3}{n+3}\big[h\phi^+_{n-2+\frac{1+(-1)^n}{2}}(\sqrt{1-8\sqrt{h}})+\phi^+_{n+2+\frac{1+(-1)^n}{2}}(\sqrt{1-8\sqrt{h}})\big].
\end{cases}\end{eqnarray}
It follows from the induction hypothesis that
\begin{eqnarray}\begin{aligned}
\sum\limits_{i+j=0}^{n+1}\xi_{i,j}I_{i,j}(h)=&\sum\limits_{i+j=0}^{n}\xi_{i,j}I_{i,j}(h)+\sum\limits_{i+j=n+1}\xi_{i,j}I_{i,j}(h)\\
=&\Big(\sum\limits_{i=0}^{[\frac{n}{4}]}\tilde{\alpha}^+_ih^i\Big) I_{0,0}(h)+\Big(\sum\limits_{i=0}^{[\frac{n-1}{4}]}\tilde{\beta}^+_ih^i\Big)I_{2,0}(h)\\
&+\Big(\sum\limits_{i=0}^{[\frac{n-3}{4}]}\tilde{\gamma}^+_ih^i\Big)I_{3,0}(h)+\sum\limits_{i+j=n+1}\xi_{i,j}I_{i,j}(h)
\\&+h\phi^+_{n-2+\frac{1+(-1)^n}{2}}(\sqrt{1-8\sqrt{h}})+\phi^+_{n+2+\frac{1+(-1)^n}{2}}(\sqrt{1-8\sqrt{h}}),
\end{aligned}\end{eqnarray}
where $\tilde{\alpha}^+_i$, $\tilde{\beta}^+_i$ and $\tilde{\gamma}^+_i$ are constants.
By substituting  (3.19) and (3.21) into (3.22) and taking into account the variation in the degrees of the coefficient polynomials associated with the integrals $I_{0,0}(h)$, $I_{2,0}(h)$ and $I_{3,0}(h)$, one derives that
\begin{eqnarray}\begin{aligned}
\sum\limits_{i+j=0}^{n+1}\xi_{i,j}I_{i,j}(h)=&\Big(\sum\limits_{i=0}^{[\frac{n+1}{4}]}\alpha^+_ih^i\Big) I_{0,0}(h)+\Big(\sum\limits_{i=0}^{[\frac{n}{4}]}\beta^+_ih^i\Big)I_{2,0}(h)+\Big(\sum\limits_{i=0}^{[\frac{n-2}{4}]}\gamma^+_ih^i\Big)I_{3,0}(h)
\\&+h\phi^+_{n-2+\frac{1+(-1)^n}{2}}(\sqrt{1-8\sqrt{h}})+\phi^+_{n+2+\frac{1+(-1)^n}{2}}(\sqrt{1-8\sqrt{h}}),
\end{aligned}\end{eqnarray}
where $\alpha^+_i$, $\beta^+_i$ and $\gamma^+_i$ are constants. This implies the validity of (3.10).

Using the same argument as (3.10), one can obtain that
\begin{eqnarray}\begin{aligned}
\sum\limits_{i+j=0}^{n}\eta_{i,j}J_{i,j}(h)=&\Big(\sum\limits_{i=0}^{[\frac{n}{4}]}\alpha^-_ih^i\Big) J_{0,0}(h)+\Big(\sum\limits_{i=0}^{[\frac{n-1}{4}]}\beta^-_ih^i\Big)J_{2,0}(h)+\Big(\sum\limits_{i=0}^{[\frac{n-3}{4}]}\gamma^-_ih^i\Big)J_{3,0}(h)
\\&+h\phi^-_{n-3+\frac{1-(-1)^n}{2}}(\sqrt{1-8\sqrt{h}})+\phi^-_{n+1+\frac{1-(-1)^n}{2}}(\sqrt{1-8\sqrt{h}}),
\end{aligned}\end{eqnarray}
where $\alpha^-_i$, $\beta^-_i$ and $\gamma^-_i$ are constants, and $\phi^-_{l}(\cdot)$ is an odd polynomials of degree no more than $l$. Plugging (3.10) and (3.24) back into (3.9) leads to (3.8). This completes the proof.\quad $\lozenge$

\vskip 0.2 true cm

From Lemma 3.2, it can be seen that in order to obtain the expression for $I(h)$, one must compute the generator integrals $I_{0,0}(h)$, $I_{2,0}(h)$, $I_{3,0}(h)$, $J_{0,0}(h)$, $J_{2,0}(h)$ and $J_{3,0}(h)$. To this end, we present the following lemma, see formula 310.05  in \cite{BF}, which play an important role in computing these integrals.
\vskip 0.2 true cm

\noindent
{\bf Lemma 3.3}\, {\it For $n\in\mathds{N}$, let $$T_n(z)=\int\frac{z^{2n}}{\sqrt{(1-z^2)(1-k^2z^2)}}dz,\ k\in(-1,1).$$ Then the following iterative formula holds:}
\begin{eqnarray}\begin{aligned}
T_{n}(z)=&\frac{z^{2n-3}}{(2n-1)k^2}\sqrt{(1-z^2)(1-k^2z^2)}+\frac{2(n-1)(k^2+1)}{(2n-1)k^2}T_{n-1}(z)\\&-\frac{2n-3}{(2n-1)k^2}T_{n-2}(z),\ n\geq2.
\end{aligned}\end{eqnarray}

For the convenience of subsequent discussion, we first present the incomplete elliptic integrals of the first, second, and third kinds in the Legendre's normal form as follows
\begin{eqnarray}\begin{aligned}
&F(z,k)=\int\frac{1}{\sqrt{(1-z^2)(1-k^2z^2)}}dz,\\
&E(z,k)=\int\frac{\sqrt{1-k^2z^2}}{\sqrt{1-z^2}}dz,\\
&\Pi(z,\mu,k)=\int\frac{1}{(1-\mu z^2)\sqrt{(1-z^2)(1-k^2z^2)}}dz,
\end{aligned}\end{eqnarray}
where $k$ and $\mu$ have the same meanings as those in the complete elliptic integrals. The following two lemmas demonstrate that $I_{0,0}(h)$, $I_{2,0}(h)$, $I_{3,0}(h)$, $J_{0,0}(h)$, $J_{2,0}(h)$ and $J_{3,0}(h)$  can be represented as a linear combination of complete elliptic integrals.

\vskip 0.2 true cm

\noindent
{\bf Lemma 3.4}\, {\it  The following equalities hold:
\begin{eqnarray}\begin{aligned}
I_{0,0}(h)=&\frac{\sqrt{1-64h}}{32a}\big[4a^2E(k)-a^2k^2(a^2-4)K(k)\\
&+(a^2-4)(a^2k^2-4)\Pi(\frac{4}{a^2},k)\big]+\frac14\sqrt{1-8\sqrt{h}},\\
I_{2,0}(h)=&\frac{\sqrt{1-64h}}{24ak^2}\big[(k^2-1)K(k)+(k^2+1)E(k)\big]+\phi^+_3\big(\sqrt{1-8\sqrt{h}}\big),\\
I_{3,0}(h)=&\frac{\sqrt{1-64h}}{240k^4a^3}\big[(k^2-1)(5k^2a^2-4k^2+8)K(k)\\&+(5k^4a^2+5k^2a^2+8k^4-8k^2+8)E(k)\big]+\phi^+_5\big(\sqrt{1-8\sqrt{h}}\big),
\end{aligned}\end{eqnarray}
where \begin{eqnarray}a=\frac{2}{\sqrt{1-8\sqrt{h}}},\ \, k=\sqrt{\frac{1-8\sqrt{h}}{1+8\sqrt{h}}}.\end{eqnarray}}

\vskip 0.2 true cm

\noindent
{\bf Proof}\, Let us commence the proof by calculating the following indefinite integral
\begin{eqnarray}\int x^{-1}\sqrt{x^4-2x^3+x^2-4h}dx.\end{eqnarray}
A simple substitution $t=x-\frac12$ transforms integral (3.29) into
\begin{eqnarray}\begin{aligned}
&\int x^{-1}\sqrt{x^4-2x^3+x^2-4h}dx\\&\qquad\quad%=\frac{\sqrt{1-64h}}{4}\int(t+\frac12)^{-1}\sqrt{(1-a^2t^2)(1-b^2t^2)}dt\\
=\frac{\sqrt{1-64h}}{2}\int\frac{1}{1-4t^2}\sqrt{(1-a^2t^2)(1-b^2t^2)}dt\\
&\qquad\quad\ \ \ \, -\sqrt{1-64h}\int\frac{t}{1-4t^2}\sqrt{(1-a^2t^2)(1-b^2t^2)}dt,
\end{aligned}\end{eqnarray}
where $a$ is as defined in (3.28) and $b=\frac{2}{\sqrt{1+8\sqrt{h}}}.$
It is noted that the integrand of the second integral on the right hand side of (3.30) is an odd function of $t$. By symmetry, the definite integral of this function with respect to $t$ over the interval $\big(-\frac12\sqrt{1-8\sqrt{h}},\frac12\sqrt{1-8\sqrt{h}}\big)$ is zero. Thus, this integral has no effect on the value of $I_{0,0}(h)$. For simplicity, we will no longer consider this integral in our subsequent transformations.
Then, changing the variable of integration in (3.30) by letting $t=\frac z a$, one obtains that
\begin{eqnarray}\begin{aligned}
%\int x^{-1}\sqrt{x^4-2x^3+x^2-4h}dx
\int\frac{1}{1-4t^2}\sqrt{(1-a^2t^2)(1-b^2t^2)}dt=&-\frac{ak^2}{4}T_1(z)-\frac{a}{16}(a^2k^2-4k^2-4)F(z,k)\\
&+\frac{(a^2-4)(a^2k^2-4)}{16a}\Pi(z,\frac{4}{a^2},k).
\end{aligned}\end{eqnarray}
Some routine calculations give rise to
\begin{eqnarray}\begin{aligned}
T_1(z)=k^{-2}F(z,k)-k^{-2}E(z,k).
\end{aligned}\end{eqnarray}
Substituting (3.32) into (3.31) yields
\begin{eqnarray}\begin{aligned}
&\frac{\sqrt{1-64h}}{2}\int\frac{1}{1-4t^2}\sqrt{(1-a^2t^2)(1-b^2t^2)}dt\\ &\qquad\qquad=\frac{a\sqrt{1-64h}}{8}E(z,k)-\frac{ak^2\sqrt{1-64h}}{32}(a^2-4)F(z,k)\\
&\qquad\qquad+\frac{(a^2-4)(a^2k^2-4)\sqrt{1-64h}}{32a}\Pi(z,\frac{4}{a^2},k)\\
&\qquad\qquad:= \Phi_1(z).
\end{aligned}\end{eqnarray}

With the above extensive preparations in place, we are now in a position  to compute $I_{0,0}(h)$. It is easy to check that as $x$ varies from $\frac12-\frac12\sqrt{1-8\sqrt{h}}$ to
$\frac12+\frac12\sqrt{1-8\sqrt{h}}$, $z$ changes from -1 to 1. Thus, one obtains through a tedious calculation that
\begin{eqnarray}\begin{aligned}
I_{0,0}(h)
=&\frac12\int_{x_0}^{x_1}x^{-1}\sqrt{x^4-2x^3+x^2-4h}dx
+\frac12\int_{x_0}^{x_1}(1-x)dx\\
=&\frac12\big(\Phi_1(1)-\Phi_1(-1)\big)+\frac14\sqrt{1-8\sqrt{h}}\\
=&\frac{a\sqrt{1-64h}}{8}E(k)-\frac{ak^2\sqrt{1-64h}}{32}(a^2-4)K(k)\\
&+\frac{(a^2-4)(a^2k^2-4)\sqrt{1-64h}}{32a}\Pi(\frac{4}{a^2},k)+\frac14\sqrt{1-8\sqrt{h}},
\end{aligned}\end{eqnarray}
where $x_0$ and $x_1$ are defined as in (3.3).

Next, we turn to calculating integrals $I_{2,0}(h)$ and $I_{3,0}(h)$. Analogous to the calculation of the indefinite integral (3.30), using the same variable substitutions as that employed there, one obtains that
\begin{eqnarray}\begin{aligned}\int x\sqrt{x^4-2x^3+x^2-4h}dx\xlongequal{t=x-\frac12}&\frac{\sqrt{1-64h}}{8}\int\sqrt{(1-a^2t^2)(1-b^2t^2)}dt\\&+
\frac{\sqrt{1-64h}}{4}\int t\sqrt{(1-a^2t^2)(1-b^2t^2)}dt,\end{aligned}\end{eqnarray}and
\begin{eqnarray}\begin{aligned}
\int x^2\sqrt{x^4-2x^3+x^2-4h}dx\xlongequal{t=x-\frac12}&\frac{\sqrt{1-64h}}{16}\int\sqrt{(1-a^2t^2)(1-b^2t^2)}dt\\&+
\frac{\sqrt{1-64h}}{4}\int t^2\sqrt{(1-a^2t^2)(1-b^2t^2)}dt.
\end{aligned}\end{eqnarray}
Similarly to the second integral on the right hand side of (3.30), we will not consider the second integral on the right hand side of (3.35) in what follows, since the definite integral of its integrand with respect to $t$ over the interval $\big(-\frac12\sqrt{1-8\sqrt{h}},\frac12\sqrt{1-8\sqrt{h}}\big)$ is also zero.

A straightforward calculation using (3.25) and (3.26) yields
\begin{eqnarray}\begin{aligned}
&\frac{\sqrt{1-64h}}{8}\int\sqrt{(1-a^2t^2)(1-b^2t^2)}dt\\&\quad\xlongequal{t=\frac za}\frac{\sqrt{1-64h}}{8a}\big[k^2T_2(z)-(k^2+1)T_1(z)+F(z,k)\big]\\
&\quad=\frac{\sqrt{1-64h}}{24a}\Big[z\sqrt{(1-z^2)(1-k^2z^2)}+\frac{k^2+1}{k^2}E(z,k)+\frac{k^2-1}{k^2}F(z,k)\Big]\\
&\quad:=\Phi_2(z),\end{aligned}\end{eqnarray}
and
\begin{eqnarray}\begin{aligned}
&\int t^2\sqrt{(1-a^2t^2)(1-b^2t^2)}dt\xlongequal{t=\frac za}\frac{k^2}{a^3}T_3(z)-\frac{k^2+1}{a^3}T_2(z)+\frac{1}{a^3}T_1(z)\\
&\qquad=\frac{(3k^2z^2-k^2-1)z}{15a^3k^2}\sqrt{(1-z^2)(1-k^2z^2)}\\&
\qquad\quad+\frac{2(k^4-k^2+1)}{15a^3k^4}E(z,k)-\frac{k^4-3k^2+2}{15a^3k^4}F(z,k).
\end{aligned}\end{eqnarray}
Substituting (3.37) and (3.38) into  (3.36) leads to
\begin{eqnarray}\begin{aligned}
 &\int x^2\sqrt{x^4-2x^3+x^2-4h}dx=\frac{\sqrt{1-64h}}{240}\Big[\frac{5k^4a^2+5k^2a^2+8k^4-8k^2+8}{a^3k^4}E(z,k)\\&
\qquad\quad+\frac{(12k^2z^2+5k^2a^2-4k^2-4)z}{a^3k^2}\sqrt{(1-z^2)(1-k^2z^2)}
\\&\qquad\quad+\frac{(k^2-1)(5k^2a^2-4k^2+8)}{a^3k^4}F(z,k)\Big]\\&
\qquad\quad:=\Phi_3(z).
\end{aligned}\end{eqnarray}

Similar to the calculation of $I_{0,0}(h)$, one can obtain that
\begin{eqnarray}\begin{aligned}
I_{2,0}(h)%=&-\int_{\Gamma^+_h}x^2ydx+\phi^+_3\big(\sqrt{1-8\sqrt{h}}\big)\\&
=&\frac12\int_{x_0}^{x_1}x\sqrt{x^4-2x^3+x^2-4h}dx\\
&+\frac12\int_{x_0}^{x_1}(x^2-x^3)dx
+\phi^+_3\big(\sqrt{1-8\sqrt{h}}\big)\\
=&\frac12\big(\Phi_2(1)-\Phi_2(-1)\big)+\frac{1}{16}\sqrt{1-8\sqrt{h}}\\&-\frac{1}{48}(1-8\sqrt{h})^\frac32+\phi^+_3\big(\sqrt{1-8\sqrt{h}}\big)\\
=&\frac{\sqrt{1-64h}}{24ak^2}\big[(k^2-1)K(k)+(k^2+1)E(k)\big]+\phi^+_3\big(\sqrt{1-8\sqrt{h}}\big),
\end{aligned}\end{eqnarray}and
\begin{eqnarray}\begin{aligned}
I_{3,0}(h)%=&-\int_{\Gamma^+_h}x^3ydx+\phi^+_5\big(\sqrt{1-8\sqrt{h}}\big)\\
=&\frac12\int_{x_0}^{x_1}x^2\sqrt{x^4-2x^3+x^2-4h}dx\\
&+\frac12\int_{x_0}^{x_1}(x^3-x^4)dx
+\phi^+_5\big(\sqrt{1-8\sqrt{h}}\big)\\
=&\frac12\big(\Phi_3(1)-\Phi_3(-1)\big)+\frac{1}{32}\sqrt{1-8\sqrt{h}}\\&-\frac{1}{160}(1-8\sqrt{h})^\frac52+\phi^+_5\big(\sqrt{1-8\sqrt{h}}\big)\\
=&\frac{\sqrt{1-64h}}{240k^4a^3}\big[(k^2-1)(5k^2a^2-4k^2+8)K(k)\\&+(5k^4a^2+5k^2a^2+8k^4-8k^2+8)E(k)\big]+\phi^+_5\big(\sqrt{1-8\sqrt{h}}\big).
\end{aligned}\end{eqnarray}
It should be noted that the final $\phi^+_3$ in (3.40) is actually equal to $$\frac{1}{16}\sqrt{1-8\sqrt{h}}-\frac{1}{48}(1-8\sqrt{h})^\frac32+\phi^+_3\big(\sqrt{1-8\sqrt{h}}\big).$$
In order to avoid introducing an excessive number of symbols, we still denote it by $\phi^+_3$. Similar treatment is also adopted in the calculation of $I_{3,0}(h)$.
The proof is complete. $\square$

\vskip 0.2 true cm

\noindent
{\bf Lemma 3.5}\, {\it  The following equalities hold:
\begin{eqnarray}\begin{aligned}
J_{0,0}(h)=&\frac{\sqrt{1-64h}}{32a}\big[4a^2E(k)-a^2k^2(a^2-4)K(k)\\
&+(a^2-4)(a^2k^2-4)\Pi(\frac{4}{a^2},k)\big]-\frac14\sqrt{1-8\sqrt{h}},\\
J_{2,0}(h)=&\frac{\sqrt{1-64h}}{24ak^2}\big[(k^2-1)K(k)+(k^2+1)E(k)\big]-\phi^+_3\big(\sqrt{1-8\sqrt{h}}\big),\\
J_{3,0}(h)=&\frac{\sqrt{1-64h}}{240k^4a^3}\big[(k^2-1)(5k^2a^2-4k^2+8)K(k)\\&+(5k^4a^2+5k^2a^2+8k^4-8k^2+8)E(k)\big]-\phi^+_5\big(\sqrt{1-8\sqrt{h}}\big),
\end{aligned}\end{eqnarray}
where $a$ and $k$ are given by (3.28)}
\vskip 0.2 true cm

\noindent
{\bf Proof} A straightforward calculation implies
$$\begin{aligned}
I_{i,j}(h)+J_{i,j}(h)=&\frac{1}{i+1}\int_{\Gamma^+_h}x^{i+1}y^jdy+\frac{1}{i+1}\int_{\Gamma^-_h}x^{i+1}y^jdy\\
=&\frac{1}{(j+1)2^{j+1}}\int_{x_0}^{x_1}x^i\big(1-x+x^{-1}\sqrt{x^4-2x^3+x^2-4h}\big)^{j+1}dx\\
&-\frac{1}{(j+1)2^{j+1}}\int_{x_0}^{x_1}x^i\big(1-x-x^{-1}\sqrt{x^4-2x^3+x^2-4h}\big)^{j+1}dx,\end{aligned}$$
where $x_0$ and $x_1$ are defined as in (3.3).
Thus, one has
\begin{eqnarray}\begin{aligned}
&I_{0,0}(h)+J_{0,0}(h)=\int_{x_0}^{x_1}x^{-1}\sqrt{x^4-2x^3+x^2-4h}dx,\\
&I_{2,0}(h)+J_{2,0}(h)=\int_{x_0}^{x_1}x\sqrt{x^4-2x^3+x^2-4h}dx,\\
&I_{3,0}(h)+J_{3,0}(h)=\int_{x_0}^{x_1}x^{2}\sqrt{x^4-2x^3+x^2-4h}dx.
\end{aligned}\end{eqnarray}
Therefore, (3.34), (3.40) and (3.41) together with (3.43) yield the desired result immediately. The proof is complete. $\square$ \vskip 0.2 true cm

As can be seen from the preceding analysis, the parameters $k$ and $a$ in $I_{0,0}(h)$, $I_{2,0}(h)$, $I_{3,0}(h)$, $J_{0,0}(h)$, $J_{2,0}(h)$ and $J_{3,0}(h)$ are all functions of $h$.
To unify the variables, we perform the variable transformation (3.6), thereby reducing them to
\begin{eqnarray}\begin{aligned}
I_{0,0}(h)=&\frac{1}{4\sqrt{2-u^2}}\Big[(2-u^2)E\big(\frac{u}{\sqrt{2-u^2}}\big)+(u^2-1)K\big(\frac{u}{\sqrt{2-u^2}}\big)\\
&+(u^2-1)^2\Pi\big(u^2,\frac{u}{\sqrt{2-u^2}}\big)\Big]+\frac{u}{4},\\
I_{2,0}(h)=&\frac{\sqrt{2-u^2}}{24}\Big[(u^2-1)K\big(\frac{u}{\sqrt{2-u^2}}\big)+E\big(\frac{u}{\sqrt{2-u^2}}\big)\Big]+\phi_3(u),\\
I_{3,0}(h)=&\frac{\sqrt{2-u^2}}{80}\Big[(u^4-2u^2+3)E\big(\frac{u}{\sqrt{2-u^2}}\big)\\&-(u^4-4u^2+3)K\big(\frac{u}{\sqrt{2-u^2}}\big)\Big]+\phi_5(u),\\
J_{0,0}(h)=&\frac{1}{4\sqrt{2-u^2}}\Big[(2-u^2)E\big(\frac{u}{\sqrt{2-u^2}}\big)+(u^2-1)K\big(\frac{u}{\sqrt{2-u^2}}\big)\\
&+(u^2-1)^2\Pi\big(u^2,\frac{u}{\sqrt{2-u^2}}\big)\Big]-\frac{u}{4},\\
J_{2,0}(h)=&\frac{\sqrt{2-u^2}}{24}\Big[(u^2-1)K\big(\frac{u}{\sqrt{2-u^2}}\big)+E\big(\frac{u}{\sqrt{2-u^2}}\big)\Big]-\phi_3(u),\\
J_{3,0}(h)=&\frac{\sqrt{2-u^2}}{80}\Big[(u^4-2u^2+3)E\big(\frac{u}{\sqrt{2-u^2}}\big)\\&-(u^4-4u^2+3)K\big(\frac{u}{\sqrt{2-u^2}}\big)\Big]-\phi_5(u).
\end{aligned}\end{eqnarray}
%To unify the variables, we perform the variable transformation $u=\sqrt{1-8\sqrt{h}}, u\in(0,1)$, thereby reducing them to
%\begin{eqnarray}\begin{aligned}
%I_{0,0}(h)=&\frac{1}{2\sqrt{2-u}}\Big[(2-u)E\big(\sqrt{\frac{u}{2-u}}\big)+(u-1)K\big(\sqrt{\frac{u}{2-u}}\big)\\ &+(u-1)^2\Pi\big(u,\sqrt{\frac{u}{2-u}}\big)\Big],\\
%I_{2,0}(h)=&\frac{\sqrt{2-u}}{12}\Big[(u-1)K\big(\sqrt{\frac{u}{2-u}}\big)+E\big(\sqrt{\frac{u}{2-u}}\big)\Big],\\
%I_{3,0}(h)=&\frac{\sqrt{2-u}}{40}\Big[(u^2-2u+3)E\big(\sqrt{\frac{u}{2-u}}\big)\\&-(u-1)(u-3)K\big(\sqrt{\frac{u}{2-u}}\big)\Big].
%\end{aligned}\end{eqnarray}

From Lemma 3.2 and (3.44), one can obtain that the first order Melnikov function $I(h)$ can be represented as a linear combination of complete elliptic integrals of first, second and third kinds and an odd polynomial  as described below.
\vskip 0.2 true cm

\noindent
{\bf Proposition 3.1}\, {\it For $u\in(0,1)$, then the first order Melnikov function $I(h)$ can be written as
\begin{eqnarray}\begin{aligned}
I(h)=\begin{cases}&\frac{1}{\sqrt{2-u^2}}\Big[P_{2[\frac{n-3}{4}]+4}(u^2)K\big(\frac{u}{\sqrt{2-u^2}}\big)+P_{2[\frac{n-3}{4}]+3}(u^2)E\big(\frac{u}{\sqrt{2-u^2}}\big)\\&\qquad+
P_{2[\frac{n}{4}]+2}(u^2)\Pi\big(u^2,\frac{u}{\sqrt{2-u^2}}\big)\Big]+\phi_{n+2}(u),\quad n\ \textup{odd},\\
&\frac{1}{\sqrt{2-u^2}}\Big[P_{2[\frac{n-1}{4}]+3}(u^2)K\big(\frac{u}{\sqrt{2-u^2}}\big)+P_{2[\frac{n-3}{4}]+3}(u^2)E\big(\frac{u}{\sqrt{2-u^2}}\big)\\&\qquad+
P_{2[\frac{n}{4}]+2}(u^2)\Pi\big(u^2,\frac{u}{\sqrt{2-u^2}}\big)\Big]+\phi_{n+1}(u), \quad n\ \textup{even},\\
\end{cases}
\end{aligned}\end{eqnarray}
where the relationship between the variables $u$ and $h$ is defined by (3.6), $P_l(u^2)$ is a polynomial in $u^2$ of degree no more than $l$ and $\phi_{l}(u)$ is an odd  polynomial in $u$ of degree no more than $l$.  }
\vskip 0.2 true cm

\noindent
{\bf Proof}\, Substituting (3.44) into (3.8) and carrying out a lengthy but routine calculation leads to the validity of (3.45).
 The proof is complete. $\Box$\vskip 0.2 true cm

 \subsection{Proof of Theorem 1.5}

In order to estimate the number of zeros of $I(h)$ on the interval $(0,\frac{1}{64})$, we will use a result introduced by Gasull, Li, Llibre and Zhang in \cite{GLLZ}, which we state in Lemma 3.6 below.
\vskip 0.2 true cm

\noindent
{\bf Lemma 3.6}  {\it Let $f_1(x)$ and $f_2(x)$ be analytic functions on an interval $\Sigma\subset \mathds{R}$. Then
$$\begin{aligned}
&\#\{x\in\Sigma|f_1(x)+f_2(x)=0\}\\&\quad\qquad\leq\#\{x\in\Sigma|f_1(x)=0\}+\#\{x\in\Sigma|f_1(x)f'_2(x)-f'_1(x)f_2(x)=0\}+1.
\end{aligned}$$}

\vskip 0.2 true cm

\noindent
{\bf Proof of Theorem 1.5}\,  Since $\phi_{n+2}(u)$ and $\phi_{n+1}(u)$ in (3.45)  are odd polynomials of $u$, we can convert it into even polynomials of $u$ by taking their first derivatives.
Here we only provide a proof for the case where $n$ is odd, and the proof for the even case can be carried out analogously. In what follows, we always denote by $P_l(\cdot)$ a polynomial of degree at most $l$.

When $n$ is an odd number, a straightforward calculation leads to
\begin{eqnarray*}
\begin{aligned}
\frac{dI(u)}{du}=&\frac{1}{(2-u^2)^{\frac32}(u-u^3)}\Big[P_{2[\frac{n+1}{4}]+4}(u^2)K\big(\frac{u}{\sqrt{2-u^2}}\big)\\&+P_{2[\frac{n+1}{4}]+3}(u^2)E\big(\frac{u}{\sqrt{2-u^2}}\big)
+P_{2[\frac{n}{4}]+4}(u^2)\Pi\big(u^2,\frac{u}{\sqrt{2-u^2}}\big)\\&+(2-u^2)^{\frac32}(u-u^3)\psi_{n+1}(u)\Big],
\end{aligned}
\end{eqnarray*}
where  $\psi_{n+1}(u)$ is an even polynomial in $u$ of degree no more than $n+1$. Introducing the change of variable $v=u^2, v\in(0,1)$, one concludes that
\begin{eqnarray}
\begin{aligned}
\frac{dI(u)}{du}%=&\frac{1}{(2-v)^{\frac32}\sqrt{v}(1-v)}\Big[P_{2[\frac{n+1}{4}]+4}(v)K\big(\sqrt{\frac{v}{2-v}}\big)+P_{2[\frac{n+1}{4}]+3}(v)E\big(\sqrt{\frac{v}{2-v}}\big)\\
%&+P_{2[\frac{n}{4}]+4}(v)\Pi\big(v,\sqrt{\frac{v}{2-v}}\big)+{(2-v)^{\frac32}\sqrt{v}(1-v)}P_{\frac{n+1}{2}}(v)\Big]
=\frac{1}{(2-v)^{\frac32}\sqrt{v}(1-v)}\big[R_0(v)+R_1(v)\big],
\end{aligned}
\end{eqnarray}
where
\begin{eqnarray*}
\begin{aligned}
R_0(v)=&{(2-v)^{\frac32}\sqrt{v}(1-v)}P_{\frac{n+1}{2}}(v),\\
R_1(v)=&P_{2[\frac{n+1}{4}]+4}(v)K\big(\sqrt{\frac{v}{2-v}}\big)+P_{2[\frac{n+1}{4}]+3}(v)E\big(\sqrt{\frac{v}{2-v}}\big)\\
&+P_{2[\frac{n}{4}]+4}(v)\Pi\big(v,\sqrt{\frac{v}{2-v}}\big).
\end{aligned}
\end{eqnarray*}
We are now in a position to estimate the number of zeros of $\frac{dI(u)}{du}$ by applying Lemma 3.6 to $R_0(v)+R_1(v)$. Some tedious computations show that
\begin{eqnarray*}
\begin{aligned}
R_0(v)R'_1(v)-R'_0(v)R_1(v)=&\sqrt{\frac{2-v}{v}}\Big[P_{2[\frac{n+1}{4}]+\frac{n+1}{2}+6}(v)K\big(\sqrt{\frac{v}{2-v}}\big)\\
&+P_{2[\frac{n+1}{4}]+\frac{n+1}{2}+5}(v)E\big(\sqrt{\frac{v}{2-v}}\big)\\
&+P_{2[\frac{n}{4}]+\frac{n+1}{2}+6}(v)\Pi\big(v,\sqrt{\frac{v}{2-v}}\big)\Big].
\end{aligned}
\end{eqnarray*}
Changing the variable of function by letting $$k=\sqrt{\frac{v}{2-v}},\ k\in (0,1),$$ one finds that
\begin{eqnarray*}
\begin{aligned}
R_0(v)R'_1(v)-R'_0(v)R_1(v)=&\frac{R(k)}{k(1+k^2)^{2[\frac{n+1}{4}]+\frac{n+1}{2}+6}},
\end{aligned}
\end{eqnarray*}
where
\begin{eqnarray*}
\begin{aligned}
R(k)=&P_{2[\frac{n+1}{4}]+\frac{n+1}{2}+6}(k^2)K(k)+P_{2[\frac{n+1}{4}]+\frac{n+1}{2}+6}(k^2)E(k)\\
&+P_{2[\frac{n+1}{4}]+\frac{n+1}{2}+6}(k^2)\Pi\big(\frac{2k^2}{1+k^2},k\big).
\end{aligned}
\end{eqnarray*}
In view of Theorem 1.4, one obtains that
\begin{eqnarray*}\begin{aligned}
&\#\{v\in(0,1)|R_0(v)R'_1(v)-R'_0(v)R_1(v)=0\}\\&\qquad\qquad=\#\{k\in(0,1)|R(k)=0\}\leq 20\Big[\frac{n+1}{4}\Big]+5n+76.
\end{aligned}\end{eqnarray*}
Observing that $R(k)$ is an even function in $k$, it follows that
\begin{eqnarray}
\#\{v\in(0,1)|R_0(v)R'_1(v)-R'_0(v)R_1(v)=0\}\leq 10\Big[\frac{n+1}{4}\Big]+\frac{5n+1}{2}+37.
\end{eqnarray}
According to (3.46) and (3.47), and taking into account that Lemma 3.6 and Rolle's Theorem, one gets that
\begin{eqnarray}\begin{aligned}
\#\{h\in(0,\frac{1}{64})|I(h)=0\}=&\#\{k\in(0,1)|I(k)=0\}\\ \leq&\#\{v\in(0,1)|R_0(v)R'_1(v)-R'_0(v)R_1(v)=0\} \\&+\#\{v\in(0,1)|R_0(v)=0\}+2\\
 \leq &10\Big[\frac{n+1}{4}\Big]+3n+40.
\end{aligned}\end{eqnarray}

When $n$ is an even number, by the argument analogous to that used above, one finds that
\begin{eqnarray}
\#\{h\in(0,\frac{1}{64})|I(h)=0\} \leq 10\Big[\frac{n}{4}\Big]+3n+42.
\end{eqnarray}

As a direct consequence of the foregoing analysis, and in light of (3.48) and (3.49), one concludes that
\begin{eqnarray*}
\#\{h\in(0,\frac{1}{64})|I(h)=0\} \leq\Big[\frac{11n}{2}\Big]+43.
\end{eqnarray*}
The proof of Theorem 1.5 is complete. $\Box$

 \section{Conclusion}
 \setcounter{equation}{0}
\renewcommand\theequation{4.\arabic{equation}}

This paper establishes  explicit upper bounds for the number of zeros of a class of functions formed by linear combinations of complete elliptic integrals of the first, second, and third kinds, where the coefficients are polynomials. The study extends the earlier important work of Gasull Li, Llibre and Zhang \cite{GLLZ} by addressing the non trivial scenario involving all three types of complete elliptic integrals, achieved through the application of Picard-Fuchs equations and well chosen variable substitutions. The derived upper bounds are expressed in terms of the degrees of the coefficient polynomials and the parametric function $\mu(k)$, providing computable estimates for both cases where $\mu(k)$ is either a polynomial or a rational function. As an application, this paper analyzes a specific piecewise smooth perturbed Hamiltonian triangle  and derives an upper bound for the number of zeros of its first order Melnikov function.
The key approach involves analyzing the algebraic structure of the Melnikov function and reducing it to manageable combinations of elliptic integrals. This research overcomes critical technical barriers in estimating the zeros of complex integral functions, thereby enriching the existing research outcomes related to the weak Hilbert's 16th problem.

However, the lower bound for the number of zeros of the first order Melnikov function $I(h)$ is not provided in the text, mainly because it is difficult to verify the independence of the coefficients in the polynomial expressions of $K(k)$, $E(k)$ and $\Pi(\mu(k),k)$ in (3.45).
It is expected that this difficulty can be overcome in the future based on the method in \cite{YZ22}.
\vskip 0.2 true cm

\noindent
{\bf CRediT authorship contribution statement}
 \vskip 0.2 true cm

\noindent
Jihua Yang: Conceptualization,  Funding acquisition, Validation, Investigation, Methodology, Resources, Supervision, Writing-original draft,	Writing-review $\&$ editing.

 \vskip 0.2 true cm

\noindent
{\bf Acknowledgment}
 \vskip 0.2 true cm

\noindent
This work was supported by the National Natural Science Foundation of China(12161069).
\vskip 0.2 true cm

\noindent
{\bf Conflict of interest}
 \vskip 0.2 true cm

\noindent
The authors declare that they have no known competing financial interests or personal relationships that could have
appeared to influence the work reported in this paper.
\vskip 0.2 true cm

\noindent
{\bf Data Availability Statement}
 \vskip 0.2 true cm

\noindent
No data was used for the research in this article. It is pure mathematics.

\newpage

\appendix

\noindent
\section{ Appendix}
\vskip 0.2 true cm

$$\begin{aligned}
W(k)=&\frac{1}{k (k^{2}-\mu(k ))^{2}(k^{2}-1)^{2}}\Big\{(1-k^2)(k^2-\mu(k)) K(k)^{3}\\
&- \big((3k^2-5)\mu(k) +2k^2\big) K(k)^2E(k)\\
&+\big((2k^2-7)\mu(k)+4k^{4}+k^{2}\big)K(k)E(k)^2+3 \mu(k) E(k)^3\\
&+\big[(k^{2}-1)(\mu(k)^{2}+2(1-2k^2)\mu(k)+k^{2}) K(k)^{2} \\
&+\frac12\big((2-k^2)\mu(k)^2+2(2k^4-4k^2+1)\mu(k)+k^2\big)K(k)E(k)\\
&-(3 \mu(k)^2+ 2(1-5k^2)\mu(k)+4k^4+k^2)E(k)^{2}\big]\Pi(\mu(k),k)\Big\}\\
&+\frac{(1-k^2)(K(k)^2-2K(k)E(k))+E(k)^2}{2k(k^2-1)\mu(k)(k^2-\mu(k))(\mu(k)-1)}\\
&\times\big[\big(k^2-\mu(k)\big)K(k)+\mu(k)E(k)+\big(\mu(k)^2-k^2\big)\Pi(\mu(k),k)\big]\mu''(k)\\
&+\frac{(1-k^2)(K(k)^2-2K(k)E(k))+E(k)^2}{2k(k^2-1)\mu(k)^2(k^2-\mu(k))^2(\mu(k)-1)^2}\\
&\times\Big[\frac12\big(k^2-\mu(k)\big)\big(5\mu(k)^2-2(2k^2+1)\mu(k)+k^2\big)K(k)\\
&+\frac{\mu(k)}{2}\big(5\mu(k)^2-2(k^2+1)\mu(k)+k^2\big)E(k)\\&
+\frac12\big(3\mu(k)^4-10k^2\mu(k)^2+4k^2(k^2+1)\mu(k)-k^4\big)\Pi(\mu(k),k)\Big]\mu'(k)^2\\
&-\frac{\mu'(k)}{2k^2(k^2-1)^2\mu(k)^2(k^2-\mu(k))^2(\mu(k)-1)}\Big\{(k^2-1)^2(k^4-\mu(k)^2)K(k)^3\\
&+(k^2-1)\big[(7k^2-3)\mu(k)^2+(k^4-5k^2)\mu(k)-2k^6+2k^4\big]K(k)^2E(k)\\
&-(k^2-1)\big[\big(10k^4-15k^2+3\big)\mu(k)^2\\
&+\big(2k^6-8k^4+10k^2\big)\mu(k)-k^6-k^4\big]K(k)E(k)^2\\
&+\big[(k^2-1)^2\big(\mu(k)^3+5k^2\mu(k)^2-5k^2\mu(k)-k^4\big)\big(K(k)^2-2K(k)E(k)\big)\\
&-\big((3k^2-1)\mu(k)^3+(3k^4-5k^2)\mu(k)^2\\&-(7k^4-5k^2)\mu(k)+k^6+k^4\big)E(k)^2\big]\Pi(\mu(k),k)\Big\}.
\end{aligned}$$


\begin{thebibliography}{50}
\bibitem{A} V. I. Arnold, Ten problems in theory of singularities and its applications, Adv. Soviet Math. 1 (1990) 1--8.
\bibitem{AV} V. I. Arnold, Loss of stability of self-oscillations close to resonance and versal deformations of equivariant vector fields, Functional Anal.  Appl. 11 (1977) 85--92.
\bibitem{BBLN} J. L. R. Bastos, C. A. Buzzi, J. Llibre, D. D. Nonaes, Melnikov analysis in nonsmooth differential systems with nonlinear switching manifold, J. Differ. Equations 267 (2019) 3748--3767.
\bibitem{BNY} G. Binyamini, D. Novikov, S. Yakovenko, On the number of zeros of Abelian integrals A constructive solution of the infinitesimal Hilbert sixteenth problem, Invent. Math. 181 (2010) 227--289.
\bibitem{BF} P. F. Byrd, M. D. Friedman, Handbook of Elliptic Integrals for Engineers and Physicists, Springer-Verlag, Berlin, 1954.
\bibitem{C} B. C. Carlson, Special Functions of Applied Mathematics, Academic Press, New York, 1977.
\bibitem{C87} B. C. Carlson, A table of elliptic integrals of the second kind, Math. Coput. 49(180) (1987) 595--606.
\bibitem{C88} B. C. Carlson, A table of elliptic integrals of the third kind, Math. Coput. 51(183) (1988) 267--280.
\bibitem{C89} B. C. Carlson, A table of elliptic integrals: cubic cases, Math. Coput. 53(187) (1988) 327--333.
\bibitem{C91} B. C. Carlson, A table of elliptic integrals: one quadratic factor, Math. Coput. 56(193) (1991) 267--280.
\bibitem{C92} B. C. Carlson, A table of elliptic integrals: two quadratic factors, Math. Coput. 59(199) (1992) 165--180.
\bibitem{CG} B. C. Carlson, J. L. Gustafson,Asymptotic approximations for symmetric elliptic integrals, SIAM J. Math. Anls. 25(2) (1994) 288--303.
\bibitem{CZW} Y. Chang, L. Q. Zhao, Q. Y. Wang, The Poincar\'{e} bifurcation by perturbing a class of cubic Hamiltonian systems, Nonlinear Anal. RWA 82 (2025) 104246.
\bibitem{CLLZ} F. D. Chen, C. Z. Li, J. Llibre, Z. F. Zhang,  A uniform proof on the weak Hilbert's 16th problem for $n=2$, J. Differ. Equations, 221 (2006) 309--342.
\bibitem{CH}  X. Y. Chen, M. A. Han, Further study on Horozov-Iliev's method of estimating the number of limit cycles, Sci. China Math. 65 (2022) 2255--2270.
%\bibitem{CGP} B. Coll, A. Gasull, R. Prohens, Bifurcation of limit cycles from two families of ceters, Dyn. Contin. Discrete Implus, Syst. Ser. A  12 (2005) 275--287.
\bibitem{DL1} F. Dumortier, C. Z. Li, Perturbations from an elliptic Hamiltonian of degree four: (I) saddle loop and two saddle cycle, J. Differ. Equations 176 (2001) 114--157.
\bibitem{DL2} F. Dumortier, C. Z. Li, Perturbations from an elliptic Hamiltonian of degree four: (II) cuspidal loop, J. Differ. Equations 175 (2001) 209--243.
\bibitem{DL3} F. Dumortier, C. Z. Li, Perturbations from an elliptic Hamiltonian of degree four: (III) global cernter, J. Differ. Equations 188 (2003) 473--511.
\bibitem{DL4} F. Dumortier, C. Z. Li, Perturbations from an elliptic Hamiltonian of degree four: (IV) figure-eight loop, J. Differ. Equations 188 (2001) 512--554.
\bibitem{GLLZ} A. Gasull, W. G. Li, J. Llibre, Z. F. Zhang, Chebyshev property of complete elliptic integrals and its application to abelian integrals, Pacific J.  Math. 202(2) (2002) 341--361.
\bibitem{GGM16} A. Gasull, A. Geyer, F. Ma\~{n}osas, On the number of limit cycles for perturbed pendulum equations, J. Differ. Equations 261 (2016) 2141--2167.
\bibitem{GGM} A. Gasull, A. Geyer, F. Ma\~{n}osas, A Chebyshev criterion with applications, J. Differ. Equations 269 (2020) 6641--6655.
 \bibitem{GI09} L. Gavrilov, I. D. Iliev, Quadratic perturbations of quadratic codimension-four centers, J. Math. Anal. Appl. 357 (2009) 69--76.
\bibitem{GI03} L. Gavrilov, I. D. Iliev, Completle hyperelliptic integrals of the first kind and their nonoscillation,  Trans. Amer. Math. Soc. 356 (2003) 1185--1207.
\bibitem{GI} L. Gavrilov, I. D. Iliev, Two-dimensional Fuchsian systems and the Chebyshev property, J. Differ. Equations 191 (2003) 105--120.
\bibitem{G} L. Gavrilov, The infinitesimal 16th Hilbert problem in the quadratic case, Invent. Math. 143 (2001) 449--497.
\bibitem{GMV} M. Grau, F. Ma\~{n}osas, J. Villadelprat, A Chebyshev criterion for Abelian integrals, Trans. Amer. Math. Soc. 363 (2011) 109--129.
\bibitem{HS} M. A. Han, L. J. Sheng, Bifurcation of limit cycles in piecewise smooth systems via Melnikov function, J. Appl. Anal. Comput. 5 (2015) 809--815.
\bibitem{HI94} E. Horozov, I. D. Iliev, On saddle-loop bifurcation of limit cycles in perturbations of quadratic Hamiltonian systems, J. Differ. Equations 113 (1994) 84--105.
\bibitem{HI} E. Horozov, I. D. Iliev, On the number of limit cycles in perturbations of quadratic Hamiltonian systems, Proc. Lond. Math. Soc. 69 (1994) 198--224
 \bibitem{HI98} E. Horozov, I. D. Iliev, Linear estimate for the number of zeros of Abelian integrals with cubic Hamiltonians, Nonlinearity 11 (1998) 1521--1537.
%\bibitem{II} I. Iliev, Perturbations of quadratic centers, Bull. Sci. math. 22 (1998) 107--161.
 %\bibitem{I} Yu. Ilyashenko, Centennial history of Hilbert's 16th problem, Bull. Amer. Math. Soc. (N.S.) 39 (2002) 301--354.
   \bibitem{IY} Y. Ilyashenko, S. Yakovenko, Double exponential estimate for the number of zeros of complete Abelian integrals, Invent. Math. 121 (1995) 613--650.
  \bibitem{LZ} C. Z. Li, Z. H. Zhang, Remarks on 16th weak Hilbert problem for $n = 2$, Nonlinearity 15 (2002) 1975--1992.
%\bibitem{LZLZ}W. Li, Y. Zhao, C. Li, Z. Zhang, Abelian integrals for quadratic centers having almost all their orbits formed by quartics, Nonlinearity 15 (2002) 863-885.
 \bibitem{LL} T. Li, J. Llibre, Limit cycles in piecewise polynomial Hamiltonian systems allowing nonlinear switching boundaries, J. Differ. Equations 344 (2023) 405--438.
\bibitem{LH12} F. Liang, M. A. Han, Limit cycles near generalized homoclinic and double homoclinic loops in piecewise smooth systems, Chaos, Solitons, Fractals 45 (2012) 454--464.
 \bibitem{LHR} F. Liang, M. A. Han, V. G. Romanovski,  Bifurcation of limit cycles by perturbing a piecewise linear Hamiltonian
system with a homoclinic loop, Nonlinear Anal. 75 (2012) 4355--4374.
 \bibitem{LX13}  C. J. Liu, D. M. Xiao, The monotonicity of the ratio of two Abelian integrals, Trans. Amer. Math. Soc. 365 (2013) 5525--5544.
\bibitem{LX}  C. J. Liu, D. M. Xiao, The smallest upper bound on the number of zeros of Abelian integrals, J. Differ. Equations 269 (2020) 3816--3852.
\bibitem{LHL} S. S. Liu, M. A. Han, J. B. Li, Bifurcation methods of periodic orbits for piecewise smooth systems, J. Differ. Equations 275 (2021) 204--233.
\bibitem{LH} X. Liu, M. A. Han, Bifurcation of limit cycles by perturbing piecewise Hamiltonian systems, Internat. J. Bifur. Chaos Appl. Sci. Engrg 20 (2010) 1379--1390.
\bibitem{M}  Y. Markov, Limit cycles of perturbations of a class of quadratic Hamiltonian vector fields, Serdica Math. J. 22(2) (1996) 91--108.
\bibitem{MV} F. Ma\~{n}osas, J. Villadelprat, Bounding the number of zeros of certain Abelian integrals, J. Differ. Equations 251 (2011) 1656--1669.
 \bibitem{NIST} F. W. J. Olver, D. W. Lozier, R. F. Boisvert, C. W. Clark, NIST Handbook of Mathematical Functionsm, Cambridge University Press, New York, 2010.
 \bibitem{PR} C. Pessoa, R. Ribeiro, Bifurcation of limit cycles from a periodic annulus formed by a center and two saddles in piecewise linear differential system with three zones, Nonlinear Anal. RWA 80 (2024) 104171.
\bibitem{PL} L. Pontryagin, On dynamic systems close to Hamiltonian systems, Zh. Eksp. Teor. Fiz. 4 (1934) 234--238.
\bibitem{R} R. Roussarie, On the number of limit cycles which appear by perturbation of separatrix loop of planar vector fields, Bol. Soc. Bras. Math. 17(2) (1986) 67--101.
\bibitem{TH} H. H. Tian, M. A. Han, Limit cycle bifurcations of piecewise smooth near-Hamiltonian systems with a switching curve, Discrete  Cont. Dyn. Syst. B 26 (2021) 5581--5599.
\bibitem{XH} Y. Q. Xiong, M. A. Han, Limit cycles appearing from a generalized heteroclinic loop with a cusp and a nilpotent saddle, J. Differ. Equations 303 (2021) 575-607.
\bibitem{XH20} Y. Q. Xiong, M. A. Han, Limit cycle bifurcations by perturbing a class of planar quintic vector fields, J. Differ. Equations 269 (2020) 10964-10994.
\bibitem{YZ17} J. H. Yang, L. Q. Zhao, The cyclicity of period annuli for a class of cubic Hamiltonian systems with nilpotent singular points, J. Differ. Equations 263 (2017) 5554--5581.
\bibitem{YZ} J. H. Yang, L. Q. Zhao, Bounding the number of limit cycles of discontinuous differential systems by using Picard-Fuchs equations, J. Differ. Equations 264 (2018) 5734--5757.
\bibitem{YZ22} J. H. Yang, L. Q. Zhao, Bifurcation of limit cycles of a piecewise smoothHamiltonian system, Qual. Theory Dyn. Syst. 21 (2022) 142.
\bibitem{ZL} Z. Zhang, C. Li, On the number of limit cycles of a class of quadratic Hamiltonian systems under quadratic perturbations, Adv. Math. 26 (1997) 445--460 (Res. Rep. 1993, 33).
\bibitem{ZSJ} L. Q. Zhao, Z. Si, R. R. Jia,  Up to the first two order Melnikov analysis for the exact cyclicity of planar piecewise linear vector fields with nonlinear switching curve, J. Differ. Equations 416 (2025) 2255--2292.
\bibitem{ZZ} Y. L. Zhao, Z. F. Zhang, Linear estimate of the number of zeros of Abelian integrals for a kind of quartic Hamiltonians, J. Differ. Equations 155 (1999) 73--88.
\end{thebibliography}
\end{document}